\documentclass[12pt]{amsart}
\usepackage{txfonts}
\usepackage{amsfonts}
\usepackage{amsfonts}
\usepackage{amssymb,latexsym}
\usepackage{enumerate}
\usepackage{mathrsfs}
\allowdisplaybreaks[3]
\newtheorem{theorem}{Theorem}

\newtheorem{proposition}[theorem]{Proposition}
\newtheorem{lemma}[theorem]{Lemma}
\newtheorem{definition}[theorem]{Definition}
\newtheorem{corollary}[theorem]{Corollary}

\newcommand{\ul}{\underline}
\newcommand {\p}{\partial}

\numberwithin{equation}{section}
\numberwithin{theorem}{section}

\def\n{\nabla}
\newcommand{\arccot}{\mathrm{arccot}}
\usepackage{hyperref}
\hypersetup{hypertex=true,colorlinks=true,linkcolor=blue,anchorcolor=blue,citecolor=blue}

\begin{document}
	\title [A new flow in K\"ahler geometry]
	{A new flow solving the LYZ equation  in K\"ahler geometry}
	\author{Jixiang Fu}
	\address{Shanghai Center for Mathematical Sciences, Jiangwan Campus, Fudan University, Shanghai, 200438, China}
	\email{majxfu@fudan.edu.cn}
	\author{Shing-Tung Yau}
	\address{Yau Mathematical Sciences Center, Tsinghua University,
		Beijing, 100084, China}
	\email{syau@tsinghua.edu.cn}
	\author{Dekai Zhang}
	\address{Department of Mathematics, Shanghai University, Shanghai, 200444, China}
	\email{dkzhang@shu.edu.cn}

	\begin{abstract}
		We {introduced  a new  flow to the LYZ equation} on a compact K\"ahler manifold. We first show the existence of the longtime solution of the flow. We then show that under the Collins-Jacob-Yau's condition on the subsolution, the longtime solution converges to the solution of the LYZ equation,
		which was solved by Collins-Jacob-Yau \cite{cjy2020cjm} by the continuity method. Moreover, as an application of the flow, we show that on a compact K\"ahler surface, if there exists a semi-subsolution of the {LYZ} equation,  then our flow converges smoothly to a singular solution to the {LYZ}  equation away from a finite number of curves of negative self-intersection. Such a solution can be viewed as a boundary point of the moduli space of the {LYZ} solutions for a given K\"ahler metric.
	\end{abstract}
	\maketitle
	\section{Introduction}
	Let $(M, \omega)$ be a compact K\"{a}hler manifold of dimension $n$ and  $\chi$ be a real closed $(1,1)$ form. Motivated by mirror symmetry {by Leung-Yau-Zaslow \cite{lyz1999}},
	Jacob-Yau \cite{jy2017ma} initiated to study the existence of solutions of  equation:
	\begin{align}\label{LYZ1}
		\mathrm{Re}(\chi_u+\sqrt{-1}\omega)^n=
		\cot\theta_0\mathrm{Im}(\chi_u+\sqrt{-1}\omega)^n,
	\end{align}
	where
	$\theta_0$ is {determined by }the complex number $\int_M(\chi+\sqrt{-1}\omega)^n$ and $\chi_u=\chi+\sqrt{-1}\partial\bar\partial u$ for a {\sl real} smooth function $u$ on $M$.
	
	Equation (\ref{LYZ1}) is called the deformed Hermitian-Yang (dHYM) equation in the literature. 
	We now call it the LYZ equation instead of the dHYM equation.

	Let $\lambda=(\lambda_1,\ldots,\lambda_n)$ be the eigenvalues of $\chi_u$ with respect to $\omega$, and if necessary we denote $\lambda$ by $\lambda(\chi_u)$ and $\lambda_i$ by $\lambda_i(\chi_u)$ for each $1\leq i\leq n$. Then by Jacob-Yau \cite{jy2017ma} the LYZ equation has an equivalent form
	\begin{align}\label{LYZ2}
		\theta_{\omega}(\chi_u):=\sum\limits_{i=1}^n\arccot{\lambda_i}=\theta_0.
	\end{align}
	It is called supercritical if  $\theta_0\in (0, \pi)$ and hypercritical if $\theta_0\in (0, \frac{\pi}{2})$.
	
	\subsection{Previous results}
	The LYZ equation has been extensively studied by many mathematicians (\cite{ch2021im,
		ccl2020arxiv,  clt2021arxiv,
		cjy2020cjm, cpw2017cvpde, cxy2018gp, hj2021cvpde, hj2020arxiv, hy2019arxiv, hzz2020arxiv, js2020arxiv, jy2017ma, li2020arxiv, li2201arxiv, pi2019arxiv}).

	We first introduce the related results in the elliptic case.
	When $n=2$, Jacob-Yau \cite{jy2017ma} solved the equation by translating it into the complex Monge-Amp\`{e}re equation which was solved by Yau \cite{ya1978cpam}.
	When $n\geq 3$,  Collins-Jacob-Yau \cite{cjy2020cjm} solved the {{LYZ}} equation for the supercritical case
	by assuming the existence of a  subsolution $\underline u$ and an extra condition on $\ul u$. For convenience, for a smooth function $v$ on $M$ we define
	\begin{align*}
		A_0(v):=\max\limits_{M}\max\limits_{1\le j\le n}\sum\limits_{ i\neq j}\arccot{\lambda_i}(\chi_{{v}})
	\end{align*}
	and
	\begin{align*}
		B_0(v)=\max\limits_{M}\theta_{\omega}(\chi_{{v}}).
	\end{align*}
	A smooth function $\ul u$ on $M$  is called a {\sl subsolution} of LYZ equation \eqref{LYZ2} if
	$\ul u$ satisfies
	the inequality
	\begin{align}\label{A0}
		A_0(\ul u)<\theta_0.
	\end{align}
	The extra condition on $\ul u$ is
	\begin{align}\label{B01}
		B_0(\ul u)< \pi.
	\end{align}
	To be precise, Collins, Jacob and Yau proved the following
	
	\begin{theorem} [Collins-Jacob-Yau \cite{cjy2020cjm}]\label{cjy}
		Let $(M,\omega)$ be a compact K\"ahler manifold of dimension $n$ and $\chi$ a closed real $(1,1)$ form on $M$ with $\theta_0\in (0,\pi)$. Suppose there exists a subsolution
		$\ul u$ of LYZ equation (\ref{LYZ2}) in the sense of  \eqref{A0} and $\ul u$ also satisfies inequality \eqref{B01}. Then there exists a unique smooth solution of LYZ equation \eqref{LYZ2}.
	\end{theorem}
	
	Without condition \eqref{B01},
	Pingali \cite{pi2019arxiv} then solved the equation for $n=3$ and Lin \cite{li2201arxiv} solved it for $n=3, 4$.
	On the other hand,
	Lin \cite{li2020arxiv} generalized Collins-Jacob-Yau's  result to the Hermitian case $(M,\omega)$ with $\partial\bar\partial\omega=\partial\bar\partial\omega^2=0$;
	Huang-Zhang-Zhang \cite{hzz2020arxiv}  considered the solution on a compact almost Hermitian manifold
	for the hypercritical case.
	
	For the parabolic flow method, there are also several results. More precisely,
	{Jacob-Yau \cite{jy2017ma} }and {Collins-Jacob-Yau \cite{cjy2020cjm}} solved the  line bundle mean curvature flow (LBMCF)
	\begin{align}\label{LBMCF}
		\left\{ {\begin{array}{*{20}l}
				{u_t=\theta_0-\theta_{\omega}(\chi_u)} \\
				{u(0)={u_0}} \\
		\end{array}} \right.
	\end{align}
	under the assumptions:
	\begin{enumerate}
		\item[(1)]  $\theta_0\in(0,\frac{\pi}{2})$;
		\item[(2)] the existence of a subsolution $\underline{u}$ in the sense of \eqref{A0}; and
		\item[(3)]  $\theta_{\omega}(\chi_{{u_0}}) \in(0,\frac{\pi}{2})$.
	\end{enumerate}
	
	{Takahashi} \cite{ta2020ijm} proved the existence and convergence of the tangent Lagrangian phase flow (TLPF)
	\begin{align}\label{TLPF}
		\left\{ {\begin{array}{*{20}l}
				{u_t=\tan\bigl(\theta_0-\theta_{\omega}(\chi_u)\bigr)} \\
				{u(0)= u_0} \\
		\end{array}} \right.
	\end{align}
	under the same assumptions (1) and (2) of flow (\ref{LBMCF}) and the assumption:
	\begin{enumerate}
		\item[(3')]  $\theta_{\omega}(\chi_{ u_0})-\theta_0\in(-\frac{\pi}{2}, \frac{\pi}{2})$.
	\end{enumerate}
	
	Another important problem raised by Collins-Jacob-Yau \cite{cjy2020cjm} is to find a sufficient and necessary geometric condition on  the existence  of a solution of the LYZ equation. There are some important progresses made by Chen \cite{ch2021im} and Chu-Lee-Takahashi \cite{clt2021arxiv}.
	
	\subsection{Our results}
	Motivated by the concavity of $\cot \theta_{\omega}(\chi_u)$ by Chen \cite{ch2021im},
	{we consider a new flow to the LYZ equation:}
	\begin{align}\label{FYZflow}
		\left\{ {\begin{array}{*{20}l}
				{u_t=\cot\theta_{\omega}(\chi_u)-\cot\theta_0}, \\
				{u(x,0)={u_0}(x)} .
		\end{array}} \right.
	\end{align}
	Assume $u_0$ satisfies
	\begin{align}\label{B0}
		B_0(u_0)<\pi.
	\end{align}
	This condition is the same as \eqref{B01} if $u_0=\ul u$.
	
	We first prove an existence theorem of the longtime solution of flow \eqref{FYZflow}.
	
	\begin{theorem}\label{els}
		Let $(M,\omega)$ be a compact K\"{a}hler manifold and $\chi$  a closed real $(1,1 )$ form with $\theta_0\in(0, \pi)$.
		If $u_0$ satisfies inequality \eqref{B0}, then  flow \eqref{FYZflow} has  a  unique smooth longtime solution $u$.
	\end{theorem}
	
	Next we consider the convergence of longtime solution of flow \eqref{FYZflow}.
	Now we need to assume the LYZ equation has a subsolution $\ul u$  which also satisfies inequality \eqref{B01}.
	The first main result of this paper is
	
	\begin{theorem}\label{mainthm}
		Let $(M,\omega)$ be a compact K\"{a}hler manifold of dimension $n$ and $\chi$  a closed real $(1,1 )$ form with $\theta_0\in(0, \pi)$. Suppose the LYZ equation \eqref{LYZ2} has a subsolution  $\underline{u}$   in the sense of \eqref{A0} which also satisfies \eqref{B01}.
		If $u_0$ satisfies \eqref{B0}, then the longtime solution $u(x,t)$ of flow \eqref{FYZflow} converges to a smooth solution $u^{\infty}$ to the  LYZ equation:
		\begin{align*}
			\theta_{\omega}(\chi_{u^{\infty}})=\theta_0.
		\end{align*}
	\end{theorem}
	
	The extra condition (\ref{B01}) in our result is the same as the one in  Theorem \ref{cjy} which is therefore reproved. Our proof here looks like simpler than the one in \cite{cjy2020cjm}.
	On the other hand,  compared with the results in \cite{jy2017ma} and \cite{ta2020ijm}, we only need $\theta_0
	\in(0, \pi)$. Moreover, condition (3) of flow \eqref{LBMCF} or  (3') of flow \eqref{TLPF}  is stronger than condition \eqref{B01}.
	
	In addition to  the concavity of $\cot\theta_{\omega}(\chi_u)$, our flow has  two advantages: The first one is the imaginary part of the Calabi-Yau functional (see the definition in Section 2) is constant along the flow, which is the key to do the $C^0$ estimate; The second one is a subsolution $\ul u$ of equation (\ref{LYZ2}) satisfying \eqref{B01} is also a subsolution of flow (\ref{FYZflow}), which allows us to use Lemma 3 in Phong-T\^o \cite{pt2017arxiv} to do higher order estimates.
	If we can establish the similar lemma without  extra condition \eqref{B01} of $\ul u$, we then can relax condition (\ref{B01}).
	
	The second motivation of this paper is to look for applications of flow (\ref{FYZflow}).
	A smooth function $\underline{u}$ is called a {\sl semi-subsolution} of the LYZ equation  if
	\begin{equation}\label{ssl}
		A_0(\ul u)\leq \theta_0.
	\end{equation}
	In the 2-dimensional case, this condition is equivalent to
	\begin{align}\label{1semisubsolution}
		\chi_{\underline{u}}\ge\cot \theta_0\omega.
	\end{align}
	Now we restrict ourselves to this case.
	
	Assume there exists a semi-subsolution $\ul u$ of the LYZ equation and replace $\chi_{\ul u}$ by $\chi$, i.e., assume that $\ul u=0$ is a semi-subsolution. For any $B_1\in (0, \pi)$, define the set
	\begin{align}\label{HB0}
		\mathcal{H}_{B_1}=\left\{v\in C^{\infty}(M, \mathbb{R}): \theta_{\omega}(\chi_{v})\in (0, B_1)\right\}.
	\end{align}
	Then if $\theta_0\in(0, \frac{\pi}{2})$,
	the set $\mathcal H_{B_1}$ for any $B_1\in (2\theta_0, \pi )$ is non-empty, for example, $0\in \mathcal{H}_{B_1}$;
	if $\theta_0\in[\frac{\pi}{2}, \pi)$, we can prove that  the set  $\mathcal{H}_{B_1}$ for any $B_1\in (\theta_0, \pi)$ is also non-empty, see Lemma \ref{2dHnonempety}.
	
	We take a function in $\mathcal H_{B_1}$ for any $B_1\in (\theta_0,\pi)$ as $u_0$ in flow \eqref{FYZflow}. We can state the second main theorem of the paper.
	
	\begin{theorem}\label{2dtheorem}
		Let $(M, \omega)$ be a compact K\"{a}hler surface and $\chi$  a closed real $(1,1)$ form.
		Assume  $\theta_0\in (0,\pi)$ and $\chi\geq\cot\theta_0\omega$.
		Then there  exist a finite number of curves $E_i$ of negative self-intersection on $M$ such that
		the solution $u(x,t)$ of  flow \eqref{FYZflow} converges to a bounded function $u^{\infty}$ in $C^{\infty}_{loc}\left(M\setminus  \cup_i E_i\right)$ as $t$ tends to $\infty$
		with the following properties.
		\begin{enumerate}
			\item $\chi+\sqrt{-1}\partial\bar\partial u^{\infty}-\cot{B_1}\omega$ is a K\"ahler current which is smooth on $M{\setminus } \cup_i E_i$;
			\item $u^{\infty}$ satisfies the LYZ equation on $M{\setminus } \cup_i E_i$
			\begin{align}\label{2d0}
				\mathrm{Re}(\chi_{u^{\infty}}+\sqrt{-1}\omega)^2=\cot \theta_0 \mathrm{Im}(\chi_{u^{\infty}}+\sqrt{-1}\omega)^2;
			\end{align}
			\item  $\chi_{u(x,t)}$ converges to $\chi_{u^{\infty}}$ and ${u^{\infty}}$  satisfies \eqref{2d0} on $M$ in the sense of currents .
		\end{enumerate}
	\end{theorem}
	
	We note that  by assuming $\theta_0\in (0, \frac{\pi}{2})$ and $B_1\leq\frac{\pi}{2} $, Takahashi \cite{ta2021cvpde} proved the same
	convergence result of the LBMCF.
	A similar result of the J-flow was studied in Fang-Lai-Song-Weinkove \cite{flsw2014apde}.
	As done by \cite{flsw2014apde, ta2021cvpde}, we need the singular solution of the degenerate complex Monge-Amp\`ere equation \eqref{2dmongeampere2} by Eyssidieux-Guedj-Zeriahi \cite{egz2009jams},
	which will be used in the $C^0$ estimate.
	We establish a similar lemma, i.e., Lemma \ref{2dsubsolutioninequality} as Lemma 3 in \cite{pt2017arxiv}  by the semi-subsolution condition to do the gradient estimate and the second order estimate. As to the convergence of  $u_t$,  the key point is that along our flow the real part of the Calabi-Yau functional is uniformly bounded.
	In this way we can prove Theorem \ref{2dtheorem}.
	
	As an application of Theorem \ref{2dtheorem}, we  have the lower bound of the $\mathcal{J}$-functional on certain spaces, see the definition in Section 2.
	
	\begin{corollary}\label{2dcorollary1}
		Let $(M, \omega)$ be a compact K\"{a}hler surface and $\chi$  a closed real $(1,1)$ form.
		Assume that $\theta_0\in(0,\pi)$ and $\chi \ge\cot\theta_0\omega $. The $\mathcal{J}$-functional  is bounded  from below in $\mathcal{H}_{B_1}$ for any $B_1\in(\theta_0, \pi)$.
	\end{corollary}
	
	If $\theta_0\in(0, \frac{\pi}{2})$, Takahashi proved that $\mathcal{J}$ is bounded  from below in $\mathcal{H}_{\frac{\pi}{2}}$.
	
	We have mentioned that for 2 dimensional case, along our flow the real part of the Calabi-Yau functional is uniformly bounded.
	We believe that the same conclusion for the higher dimension also holds.
	Hence the real part of the Calabi-Yau functional plays the similar role as the Donaldson functional defined on the space of Hermitian metrics on
	a holomorphic vector bundle. We expect  that we can use our flow to study the moduli space of solutions of the LYZ equation
	on a compact K\"ahler manifold $(M,\omega)$.
	
	The rest of this paper is arranged as follows.
	In Section 2, we give some preliminary results on the linearized operator on our flow,
	the concavity of $\cot\theta(\lambda)$, the parabolic subsolution, and the Calabi-Yau functional.
	In Section 3, we prove Theorem \ref{els}.
	In Section 4, we prove Theorem \ref{mainthm}, including  the $C^0$ estimate, the gradient estimate and the second order estimate.
	In Section 5, we prove Theorem \ref{2dtheorem} and Corollary \ref{2dcorollary1}.
	
	\textbf{Notations:}
	In this paper a closed real $(1,1)$ form $\chi$ is fixed.
	We will use the constant $C$ in the generic sense
	which is dependent on $\omega$, $\chi$, $\underline u$, $u_0$ and $n$.
	If necessary, we will use $C_i$ as a specific constant.
	
	Notations of covariant derivatives are used. For example,  $u_{i\bar j  k}$ represents the third order covariant derivative of  function $u$, $\alpha_{i\bar j, k}$ represents the covariant derivative of   (1,1) form   $\alpha$.
	
	We use Einstein summation convention if there is no confusion.
	
	\section{Preliminary results}
	\subsection{The linearized operator}
	Note
	\begin{align}\label{060301}
		\cot\theta_{\omega}(\chi_u)=\frac{\mathrm{Re}(\chi_u+\sqrt{-1}\omega)^n}{\mathrm{Im}(\chi_u+\sqrt{-1}\omega)^n}.
	\end{align}
	We manipulate the linearized operator
	$\mathcal{{P}}$ of our flow (\ref{FYZflow}) in the following lemma.
	
	\begin{lemma}
		The operator $\mathcal P$ has the form:
		\begin{align*}
			\mathcal{{P}}(v)=v_t-\csc^2\theta_{\omega}(\chi_u)\bigl(wg^{-1}w+g\bigr)^{i\bar j}v_{i\bar j},
		\end{align*}
		where $g=(g_{i\bar j})_{n\times n}$, $w=(w_{i\bar j})_{n\times n}$ for  $w_{i\bar j}=\chi_{i\bar j}+u_{i\bar j}$, and $D^{i\bar j}:=(D^{-1})_{i\bar j}$
		for an invertible Hermitian symmetric matrix $D$.
	\end{lemma}
	
	\begin{proof}
		We only need to deal with the variation of $\cot\theta_{\omega}(\chi_u)$. Indeed, let $u(s)$ be a variation of the function $u$ and $\frac{du(s)}{ds}|_{s=0}=v$. Let $A(s):=g^{-1}w(s)+\sqrt{-1}I$ with $w(s)$ being the local matrix of $\chi_{u(s)}$. Then
		\begin{align}\label{060701}
			A(s)^{-1}=\bigl(g^{-1}w(s)-\sqrt{-1}I\bigr)\bigl((g^{-1}w(s))^2+I\bigr)^{-1}.
		\end{align}
		For simplicity, we write $A$ instead of $A(s)$.
		By (\ref{060301}) we have
		\begin{align*}
			\delta\bigl(\cot\theta_{\omega}(\chi_u)\bigr)=&\frac{\mathrm{Re}{\left( \delta\det A\right)}}{\mathrm{Im}(\det A)}-\frac{\mathrm{Re}(\det A)\mathrm{Im}{ \left(\delta \det A\right)}}{(\mathrm{Im}(\det A))^2}.
		\end{align*}
		Since  $\delta(\det A)=(\det A)\delta(\log\det A)$, if we write $\det A=a_1+\sqrt{-1}a_2$ and $\delta\left(\log\det A\right)=b_1+\sqrt{-1}b_2$, then
		\begin{align*}
			\delta \bigl(\cot\theta_{\omega}(\chi_u)\bigr)=&\frac{a_1b_1-a_2b_2}{a_2}-\frac{a_1(a_1b_2+a_2b_1)}{a_2^2}\notag\\
			=&\frac{-a_1^2-a_2^2}{a_2^2} b_2=-\csc^2\theta_{\omega}(\chi_u) b_2.
		\end{align*}
		On the other hand, by (\ref{060701}) we have
		\begin{align*}
			b_2=&\textup{Im}\,\delta(\log\det A)=-\textup{tr}\bigl((wg^{-1}w+g)^{-1}\delta w(s)|_{s=0}\bigr)\\
			=&-(wg^{-1}w+g)^{i\bar j}v_{i\bar j}.
		\end{align*}
		Hence
		\begin{align}\label{linearize3}
			\delta \bigl(\cot\theta_{\omega}(\chi_u)\bigr)=\csc^2\theta_{\omega}(\chi_u)\bigl(w g^{-1}w+g\bigr)^{i\bar j}   v_{i\bar j}.
		\end{align}
	\end{proof}
	We denote
	\begin{align}\label{Fijdef}
		F^{i\bar j}:=\textup{csc}^2\theta_{\omega}(\chi_u)\bigl(wg^{-1}w+g\bigr)^{i\bar j}
	\end{align}
	and hence
	\begin{align*}
		\mathcal P(v)=v_t-F^{i\bar j}v_{i\bar j}.
	\end{align*}
	The following lemma is useful in the gradient and second order estimates.
	
	\begin{lemma}\label{0606lemma}
		Let $u$ be a solution of  flow \eqref{FYZflow}. Then
		\begin{align}
			&u_{tp}-F^{i\bar j}w_{i\bar j, p}=0,\label{060600}\\
			& u_{tp\bar p}-F^{i\bar j}w_{i\bar j, p\bar p}\notag\\
			&= -F^{i\bar l}\bigl(wg^{-1}w+g\bigr)^{k\bar j}w_{i\bar j, p}(w_{k\bar m, \bar p}g^{r\bar m}w_{r\bar l}+w_{k\bar m}g^{r\bar m}w_{r\bar l, \bar p})\notag\\
			&+2\cot\theta_{\omega}(\chi_u)F^{i\bar j}w_{i\bar j, p}\bigl(wg^{-1}w+g\bigr)^{k\bar l}w_{k\bar l,\bar p}.\label{060601}
		\end{align}
	\end{lemma}
	
	\begin{proof}
		Similar as the proof of \eqref{linearize3}, differentiating  equation  \eqref{FYZflow} leads to (\ref{060600}) directly:
		\begin{align*}
			u_{tp}=\csc^2\theta_{\omega}(\chi_u)\bigl(w g^{-1}w+g\bigr)^{i\bar j}    w_{i\bar j, p}=F^{i\bar j} w_{i\bar j, p}.
		\end{align*}
		
		Differentiating the equation twice, we have
		\begin{align}
			u_{tp\bar p}
			=&\textup{csc}^2\theta_{\omega}(\chi_u)\bigl(wg^{-1}w+g\bigr)^{i\bar j}w_{i\bar j, p\bar p}\notag\\ &+(\csc^2\theta_{\omega}(\chi_u))_{\bar p}\bigl(wg^{-1}w+g\bigr)^{i\bar j}w_{i\bar j, p}\notag\\
			-&\textup{csc}^2\theta_{\omega}(\chi_u)\bigl(wg^{-1}w+g\bigr)^{i\bar l}\bigl(wg^{-1}w+g\bigr)^{k\bar j}w_{i\bar j, p}(wg^{-1}w+g)_{k\bar l, \bar p},\notag
		\end{align}
		where
		\begin{align*}
			(\csc^2\theta_{\omega}(\chi_u))_{\bar p}
			=2\cot\theta_{\omega}(\chi_u)(\cot\theta_{\omega}(\chi_u))_{\bar p}
			=2\cot\theta_{\omega}(\chi_u)F^{k\bar l}w_{{k\bar l, \bar p}}
		\end{align*}
		and
		\begin{align*}
			(wg^{-1}w+g)_{k\bar l, \bar p}=&\bigl(w_{k\bar m}g^{r\bar m}w_{r\bar l}+g_{k\bar l}\bigl)_{\bar p}\\
			=&w_{k\bar m, \bar p}g^{r\bar m}w_{r\bar l}+w_{k\bar m}g^{r\bar m}w_{r\bar l, \bar p}.
		\end{align*}
		Hence identity \eqref{060601} follows.
	\end{proof}
	
	\subsection{The concavity of $\cot\theta(\lambda)$ in $\Gamma_{\tau}$  for $\tau\in (0,\pi)$.}
	
	Here
	\begin{align}\label{070601}
		\theta(\lambda):=\sum_{i=1}^n\arccot\lambda_{i}\quad \textup{for}\ \lambda=(\lambda_1,\dots,\lambda_n)\in\mathbb R^n
	\end{align}
	and
	\begin{align*}
		\Gamma_{\tau}:=\{\lambda\in\mathbb{R}^n \ | \ \theta(\lambda)<\tau\}\subset\mathbb R^n
		\quad\textup{for $\tau\in(0,\pi)$}.
	\end{align*}
	We have the following useful facts.
	
	\begin{lemma}[Yuan \cite{yu2006pams}, Wang-Yuan  \cite{wy2014ajm}]\label{060905}
		If $\theta(\lambda)\leq \tau\in(0,\pi)$ for $\lambda=(\lambda_1,\dots,\lambda_n)$ with  $\lambda_1\ge \lambda_2\ge\cdots\ge  \lambda_n$, then the following inequalities holds.
		\begin{enumerate}
			\item[(1)] $ \lambda_1\ge \lambda_2\ge\cdots\ge\lambda_{n-1}\ge \cot{\frac{\tau}{2}}>0$, {and}  $\lambda_{n-1}\ge |\lambda_n|$;
			\item[(2)]
			{$\lambda_1+(n-1)\lambda_n\ge 0$.}
		\end{enumerate}
		Moreover, $\Gamma_{\tau}$ is convex for any $\tau\in(0,\pi)$.
	\end{lemma}
	
	\begin{lemma}[Chen \cite{ch2021im}]\label{convexity}
		For any $\tau\in (0, \pi)$,  the function $\cot\theta(\lambda)$ on $\Gamma_{\tau}$ is concave.
	\end{lemma}
	\begin{proof}
		For completeness, we give an elementary proof here.\\
		When $n=1$, $\cot\theta(\lambda)=\lambda_1$ is obviously concave. We now assume $n\ge 2$. By definition (\ref{070601}) we have
		\begin{align}\label{060702}
			\frac{\partial^2 \cot\theta{(\lambda)}}{\partial\lambda_i\partial\lambda_j}=&-\frac{\partial}{\partial \lambda_j}\Bigl(\csc^2\theta(\lambda)\frac{\partial\theta(\lambda)}{\partial\lambda_i}\Bigl)=
			\frac{\partial}{\partial \lambda_j}\Bigl(\frac{\csc^2\theta(\lambda)}{1+\lambda_i^2}\Bigl)\notag\\
			=&-2\csc^2\theta(\lambda)\Bigl(\frac{\lambda_i\delta_{ij}}{(1+\lambda_i^2)^2}-\frac{\cot\theta(\lambda)}{(1+\lambda_i^2)(1+\lambda_j^2)}\Bigr).
		\end{align}
		Hence the function $\cot\theta(\lambda)$ on $\Gamma_\tau$ is concave if and only if the matrix
		\begin{align*}
			\Lambda=
			\bigl( \lambda_i\delta_{ij}-\cot\theta(\lambda)   \bigr)_{n\times n}
		\end{align*}
		is positive definite.
		Without loss of generality, we assume $\lambda_1\ge \lambda_2\ge\dots\ge\lambda_n$. Since $\theta(\lambda)\in(0, \pi)$, by Lemma \ref{060905} (1), we have $\lambda_{n-1}>0$.
		
		By the definition of $\theta(\lambda)$,  for any $1\le j_1<j_2<\dots<j_k$, $1\le k\le n-1$,
		we have $\sum_{l=1}^k\arccot \lambda_{j_l}<\theta(\lambda)$. Hence
		\begin{align}\label{621}
			\mathrm{Re}\left(\prod\limits_{l=1}^k(\lambda_{j_l}
			+\sqrt{-1})\right)- \cot\theta(\lambda)\mathrm{Im}\left(\prod\limits_{l=1}^k
			(\lambda_{j_l}+\sqrt{-1})\right)>0.
		\end{align}
		Let $\sigma_i(\lambda_{j_1j_2\dots j_k})$ for $1\leq i\leq k$ be  the $i$-th elementary symmetric polynomial of $\lambda_{j_1},\lambda_{j_2},\dots,\lambda_{j_k}$. Then
		\eqref{621} can be written as
		\begin{align}\label{622}
			\sum_{i=0}^{[k/2]}(-1)^i\sigma_{k-2i}(\lambda_{j_1j_2\dots j_k})-\cot\theta(\lambda)\sum_{i=0}^{[(k-1)/2]}(-1)^i\sigma_{k-1-2i}(\lambda_{j_1j_2\dots j_k})>0.
		\end{align}
		Denote by $D_k$ the $k$-th leading principal minor of the matrix $\Lambda$.  We need to  prove  $D_k>0$ for any $1\le k \le n$.
		When $k=1$, $D_1=\lambda_1-\cot\theta(\lambda)>0$.
		When $2\le k\le n$, by direct computation, we have
		\begin{align*}
			D_k=&\sigma_k(\lambda_{12\ldots k})-\cot\theta(\lambda)\sigma_{k-1}(\lambda_{12\ldots k}).
		\end{align*}
		Hence by \eqref{622}, we have
		\begin{align*}
			D_k &>-\sum_{i=1}^{[k/2]}(-1)^i\sigma_{k-2i}(\lambda_{12\dots k})+\cot\theta(\lambda)\sum_{i=1}^{[(k-1)/2]}(-1)^i\sigma_{k-1-2i}(\lambda_{12\dots k})\notag\\
			&=:E_{k-2}(\lambda_{12\dots k})\notag
		\end{align*}
		We prove  $E_{k-2}(\lambda_{12\dots k})>0$ for any $2\le k\le n$.
		
		We use the well-known formula
		\begin{align}\label{623}
			\sigma_i(\lambda_{12\dots k})=&\sigma_i(\lambda_{2\dots k})+\lambda_1\sigma_{i-1}(\lambda_{2\dots k})
		\end{align}
		for $1\leq i\leq k-1$ to deduce that
		\begin{align}\label{624}
			E_{k-2}(\lambda_{12\dots k})=F_{k-2}(\lambda_{2\dots k})+\lambda_1 E_{k-3}(\lambda_{2\dots k}),
		\end{align}
		where
		\begin{align*}
			F_{k-2}(\lambda_{2\dots k})=&-\sum_{i=1}^{[k/2]}(-1)^i\sigma_{k-2i}(\lambda_{2\dots k})\\
			&+\cot\theta(\lambda)\sum_{i=1}^{[(k-1)/2]}(-1)^i\sigma_{k-1-2i}(\lambda_{2\dots k})\\
			=&\sum_{j=0}^{[(k-2)/2]}(-1)^j\sigma_{k-2-2j}(\lambda_{2\dots k})\\
			&-\cot\theta(\lambda)\sum_{j=0}^{[(k-3)/2]}(-1)^j\sigma_{k-3-2j}(\lambda_{2\dots k})
		\end{align*}
		and
		\begin{align*}
			E_{k-3}(\lambda_{2\dots k})=&-\sum_{i=1}^{[k/2]}(-1)^i\sigma_{k-2i-1}(\lambda_{2\dots k})\\
			&+\cot\theta(\lambda)\sum_{i=1}^{[(k-1)/2]}(-1)^i\sigma_{k-2-2i}(\lambda_{2\dots k})\\
			=&\sum_{j=0}^{[(k-2)/2]}(-1)^j\sigma_{k-3-2j}(\lambda_{2\dots k})\\
			&-\cot\theta(\lambda)\sum_{j=0}^{[(k-3)/2]}(-1)^j\sigma_{k-4-2j}(\lambda_{2\dots k}).
		\end{align*}
		
		By \eqref{623} we compute directly to get
		\begin{align}
			F_{k-2}(\lambda_{2\dots k})
			=& \mathrm{Re}\left(\prod\limits_{j=3}^{k}(\lambda_{j}+\sqrt{-1})\right)\\
			&-\cot\theta(\lambda)\mathrm{Im}\left(\prod\limits_{j=3}^{k}(\lambda_{j}+\sqrt{-1})\right)+\lambda_2F_{k-3}(\lambda_{3\dots k}).\notag
		\end{align}
		Hence
		\begin{align*}
			F_{k-2}(\lambda_{2\dots k})>\lambda_2F_{k-3}(\lambda_{3\ldots k}).
		\end{align*}
		From this we deduce that
		\begin{align*}
			F_{k-2}(\lambda_{2\dots k})>&\lambda_2\lambda_3\cdots\lambda_{k-2}F_{1}(\lambda_{(k-1) k})\\
			=&\lambda_2\lambda_3\cdots\lambda_{k-2}(\lambda_{k-1}+\lambda_k-\cot\theta(\lambda))>0.
		\end{align*}
		Combined with \eqref{624}, we have
		\begin{align*}
			E_{k-2}(\lambda_{12\dots k})>\lambda_1 E_{k-3}(\lambda_{2\dots k}).
		\end{align*}
		Hence  for any $2\leq k\leq n$ we have
		\begin{align*}
			E_{k-2}(\lambda_{12\dots k})>&\lambda_1\lambda_2\cdots\lambda_{k-3} E_{1}(\lambda_{(k-2)(k-1)k})\\
			=&\lambda_1\lambda_2\cdots\lambda_{k-3} (\lambda_{k-2}+\lambda_{k-1}+\lambda_k-\cot\theta(\lambda))
			>0.
		\end{align*}
		
		In summary, we finish the proof of the lemma.
	\end{proof}
	
	\subsection{Parabolic subsolution }
	B. Guan \cite{gu2014dmj} introduced the definition  of a subsolution of fully non-linear equations. G. Sz\'{e}kelyhidi \cite{sz2019jdg}
	gave a weaker version of a subsolution  and Collins-Jacob-Yau \cite{cjy2020cjm} used it to the LYZ equation which is equivalent to (\ref{A0}). {These two notions  are equivalent for the type 1 cones by the appendix in \cite{guan2014arxiv}.}  On the other hand, Phong-T\^o \cite{pt2017arxiv}  modified the definition in \cite{gu2014dmj} and \cite{sz2019jdg} to the parabolic case. We use their definition to our flow.
	
	\begin{definition}\label{060901} A smooth function
		$\underline{u}(x, t)$ on $M\times [0,T)$ is called a subsolution of flow (\ref{FYZflow})
		if there exists a constant $\delta>0$ such that for any $(x,t)\in M\times [0, T)$, the subset
		of  $\mathbb R^{n+1}$
		\begin{align*}
			S_{\delta}(x,t):=\{(\mu,\tau)&\in \mathbb{R}^n\times \mathbb{R}\ |\ \cot\theta\bigl(\lambda(\chi_{\underline{u}(x,t)})+\mu\bigr)-{\underline{u}}_t(x, t)+\tau\\=&\cot\theta_0, \
			\mu_i>-\delta\ \text{for  $1\le i\le n$ and} \ \tau>-\delta\}
		\end{align*}
		is uniformly bounded, i.e., it is contained in the ball $B_{K}^{{n+1}}(0)$ in $\mathbb R^{n+1}$ with radius $K$, a uniform constant.
	\end{definition}
	
	We have the following observation.
	\begin{lemma}
		If $\underline{u}$ is a subsolution of  LYZ equation \eqref{LYZ1} with $B_0(\ul u)<\pi$, then the function $\underline u(x,t)=\underline u (x)$ on $M\times [0,\infty)$ is also a subsolution of  (\ref{FYZflow}).
	\end{lemma}
	
	\begin{proof}
		We want to find a  constant $\delta$ in Definition \ref{060901}.
		If such a  $\delta$  exists, we let $(\mu,\tau)\in S_{\delta}(x,t)$ for $(x,t)\in M\times [0,\infty)$.
		Since $\mu_i>-\delta$ for each $1\leq i\leq n$, by the definition of $B_0(\underline{u})$ in (\ref{B01}) we have
		\begin{align*}
			0<\theta\bigl(\lambda(\chi_{\underline{u}(x)})+\mu\bigr)\leq \theta_{\omega}(\chi_{\underline{u}(x)})+n\delta\leq B_0(\underline{u})+n\delta.
		\end{align*}
		Hence if  $0<\delta\leq \frac{\pi-B_0(\underline{u})}{2n}$, then
		\begin{align*}
			0<\theta(\lambda(\chi_{\underline{u}(x)})+\mu)<\frac{\pi+B_0(\underline{u})}{2}<\pi,
		\end{align*}
		and by the definition of $S_\delta(x,t)$, $\tau$ is bounded from above:
		\begin{align*}
			\tau=\cot\theta_0-\cot\theta(\lambda(\chi_{\underline{u}(x)})+\mu)\le \cot\theta_0-\cot\Bigl(\frac{\pi+B_0(\underline{u})}{2}\Bigr).
		\end{align*}
		
		Since also $\mu_i>-\delta$ for each $1\leq i\leq n$,
		by subsolution condition (\ref{A0}) we  have
		\begin{align*}
			\sum\limits_{i\neq j}\arccot(\lambda_i(\chi_{\underline{u}(x)})+\mu_i)
			\le& \sum\limits_{i\neq j}\arccot\lambda_i(\chi_{\underline{u}(x)})+(n-1)\delta\\
			\le &A_0(\underline{u})+(n-1)\delta.
		\end{align*}
		If  $0<\delta\leq \frac{\theta_0-A_0(\underline{u})}{2(n+1)}$, then
		\begin{align*}
			\sum\limits_{i\neq j}\arccot(\lambda_i(\chi_{\underline{u}(x)})+\mu_i)\le \frac{\theta_0+A_0(\underline{u})}{2}.
		\end{align*}
		Since $\tau>-\delta$, by the definition of $S_\delta(x,t)$ we have for each $j$
		\begin{align*}
			\arccot(\lambda_j(\chi_{\underline{u}})+\mu_j)=&\arccot(\cot\theta_0-\tau)-\sum
			\limits_{i\neq j}\arccot(\lambda_i(\chi_{\underline{u}})+\mu_i)\\
			\ge& \arccot(\cot\theta_0+\delta)-\sum
			\limits_{i\neq j}\arccot(\lambda_i(\chi_{\underline{u}})+\mu_i)\\
			\ge& \theta_0-\delta-\frac{\theta_0+A_0(\underline{u})}{2}
			\ge\frac{n(\theta_0-A_0(\underline{u}))}{2(n+1)}>0.
		\end{align*}
		Hence  we have
		\begin{align*}
			\mu_j\le \max\limits_{ M}|\lambda(\chi_{\underline{u}(x)})|+\cot\Bigl(\frac{n(\theta_0-A_0(\underline{u}))}{2(n+1)}\Bigr).
		\end{align*}
		
		Therefore, if we choose $\delta=\min\{\frac{\pi-B_0(\underline{u})}{2n},\frac{\theta_0-A_0(\underline{u})}{2(n+2)}\}$,
		then
		for any $(x, t)\in M\times[0, \infty)$ and $(\mu, \tau)\in S_{\delta}(x, t)$,  we have
		\begin{align*}
			|\mu|+|\tau|\le K:=2n\Bigr(\delta&+\max\limits_{ M}|\lambda(\chi_{\underline{u}})|+\cot\theta_0-\cot\Bigl(\frac{\pi+B_0(\underline{u})}{2}\Bigr) \\
			&+\cot\Bigl(\frac{n(\theta_0-A_0(\underline{u}))}{2(n+1)}\Bigr)\Bigr).
		\end{align*}
	\end{proof}
	
	\subsection{The Calabi-Yau Functional}
	Recall the definition of the Calabi-Yau functional by Collins-Yau \cite{cy2018arxiv}: for any  $v\in C^{2}(M, \mathbb{R})$,
	\begin{align*}
		\textup{CY}_{\mathbb{C}}(v):=\frac{1}{n+1}\sum\limits_{i=0}^n \int_{M}{v(\chi_v+\sqrt{-1}\omega)^i\wedge(\chi+\sqrt{-1}\omega)^{n-i}}.
	\end{align*}
	The $\mathcal{J}$-functional is defined by
	\begin{align*}
		\mathcal{J}(v):=\mathrm{Im}\bigl(e^{-\sqrt{-1}\theta_0}\textup{CY}_{\mathbb{C}}(v)\bigr).
	\end{align*}
	
	Let $v(s)\in C^{2,1}(M\times[0,T], \mathbb{R})$ be a variation of the function  $v$, i.e., $v(0)=v$. The integration by parts gives
	\begin{align}\label{060302}
		\frac{d}{ds}\textup{ CY}_{\mathbb{C}}(v(s))= &\int_{M}{\frac{\partial v(s)} {\partial s}\bigl(\chi_{v(s)}+\sqrt{-1}\omega\bigr)^n},\\
		\frac{d}{ds}\mathcal{J}(v(s))=&\int_{M}{\frac{\partial v(s)} {\partial s}\mathrm{Im}(e^{-\sqrt{-1}\theta_0}\bigl(\chi_{v(s)}+\sqrt{-1}\omega\bigr)^n)}.
	\end{align}
	
	\begin{lemma}\label{IMCY}
		Let $u(x, t)$ be a solution of flow \eqref{FYZflow}. Then
		\begin{align}
			&  \mathrm{Im}\bigl(\textup{CY}_{\mathbb{C}}(u(\cdot,t))\bigr)= \mathrm{Im}\bigl(\textup{CY}_{\mathbb{C}}( u_0)\bigr)\label{CYF1},\\
			& \frac{d}{dt}\mathrm{Re}\bigl(\textup{CY}_{\mathbb{C}}(u(\cdot,t))\bigr)=\int\limits_M{{\Bigl(\frac{\partial u(t)}{\partial t}\Bigr)^2}\mathrm{Im}(\chi_u+\sqrt{-1}\omega)^n}\label{CYF2},\\
			&\frac{d}{dt}\mathcal{J}(u(\cdot,t)))\le 0\label{CYFJ}.
		\end{align}
	\end{lemma}
	
	\begin{proof}
		Denote by $u(t):=u(x, t)$ for simplicity. Then we have
		\begin{align*}
			&\frac{d}{dt}\mathrm{Im}\bigl(\textup{CY}_{\mathbb{C}}(u(t))\bigr)\\
			=&\int_{M}{\frac{\partial u(t)}{\partial t}\mathrm{Im}(\chi_{u(t)}+\sqrt{-1}\omega)^n} \\
			=&\int_{M}{\Bigl(  \frac{\mathrm{Re} (\chi_{u(t)}+\sqrt{-1} \omega)^n}{\mathrm{Im}(\chi_{u(t)}+\sqrt{-1} \omega)^n} -\cot\theta_0  \Bigr)\mathrm{Im}(\chi_{u(t)}+\sqrt{-1}\omega)^n}\\
			=&\int_{M}{  {\mathrm{Re} (\chi+\sqrt{-1} \omega)^n}} -\cot\theta_0  \int_M{\mathrm{Im}(\chi+\sqrt{-1}\omega)^n}\\
			=&0,
		\end{align*}
		where each equality is successively by (\ref{060302}), (\ref{FYZflow}) and (\ref{060301}), Stokes' theorem, and the definition of $\theta_0$. Hence  \eqref{CYF1} holds as $u(0)= u_0$.
		
		Then  we can also prove \eqref{CYF2}.
		\begin{equation*}
			\begin{aligned}
				&\frac{d}{dt}\textup{Re}(\textup{CY}_{\mathbb C}(u(t)))\\
				=&\int_M\frac{\partial u(t)}{\partial t}\textup{Re}(\chi_{u(t)}+\sqrt{-1}\omega)^n\\
				=&\int_M\frac{\partial u(t)}{\partial t}\cot\theta_{\omega}(\chi_{u(t)})\textup{Im}(\chi_{u(t)}+\sqrt{-1}\omega)^n\\
				=&\int_M\frac{\partial u(t)}{\partial t}\Bigl(\frac{\partial u(t)}{\partial t}+\cot\theta_0\Bigr)\textup{Im}(\chi_{u(t)}+\sqrt{-1}\omega)^n\\
				=&\int_M\Bigl(\frac{\partial u(t)}{\partial t}\Bigr)^2\textup{Im}(\chi_{u(t)}+\sqrt{-1}\omega)^n,
			\end{aligned}
		\end{equation*}
		where the last equality follows from  (\ref{CYF1}) .
		
		Locally
		\begin{align*}
			& \frac{\partial u(t)}{\partial t}\mathrm{Im} (
			e^{-\sqrt{-1}\theta_0}(\chi_u+\sqrt{-1}\omega)^n )\\
			=&\prod_{i=1}^{n}(1+\lambda_i^2)(\cot\theta_{\omega}(\chi_{u(t)})-\cot\theta_0)\sin(\theta_{\omega}(\chi_{u(t)})-\theta_0)\omega^{n}
			\\
			=&-\prod_{i=1}^{n}(1+\lambda_i^2)\frac{\sin^2(\theta_{\omega}(\chi_{u(t)})-\theta_0)}{\sin\theta_{\omega}(\chi_{u(t)})\sin\theta_0}\omega^{n}\le 0,
		\end{align*}
		where the last inequality follows from $\theta_{\omega}(\chi_{u(t)})\in (0, \pi)$ by \eqref{060902}. Hence $\mathcal{J}$ is decreasing and \eqref{CYFJ} follows.
	\end{proof}
	
	Next we prove that along our flow the real part of the Calabi-Yau functional can be controlled by $|u|_{L^{\infty}}$ without  the subsolution condition.
	\begin{proposition}\label{recy}
		Let $u(x, t)$ be a solution of  flow \eqref{FYZflow} with the initial data satisfying $\eqref{B0}$. Then there exists a uniform constant $C$ such that
		\begin{align}
			\mathrm{Re}\left(\textup{ CY}_{\mathbb{C}}(u)\right)\le C |u|_{L^{\infty}}.
		\end{align}
	\end{proposition}
	\begin{proof}
		By the definition of  the Calabi-Yau functional, we only need to prove that for any $0\le k,l\le n$ with $0\le k+l\le n$
		\begin{align}\label{ReCY1}
			\left|\int_M{u\ \chi_u^{k}\wedge  \chi^{l}\wedge\omega^{n-k-l} }\right| \le C |u|_{L^{\infty}}.
		\end{align}
		We prove the above estimates by inductive argument on $k$. When  $k=0$, it obviously holds. Now assume inequality \eqref{ReCY1} holds for $k\le m$ with
		$0\le k+l\le n$. We prove inequality \eqref{ReCY1} holds for $k= m+1$.
		Indeed, since along the flow by \eqref{060902} $\chi_{u}\ge -
		\cot B_0(u_0)\omega$, there exists a constant $C_0>0$ such that  $\chi_u+C_0\omega>0$  and $\chi+C_0\omega>0$. We write
		\begin{equation*}
			\begin{aligned}
				\int_M &u\chi_u^{m+1}\wedge\chi^l\wedge\omega^{n-m-l-1}\notag\\
				=&\int_M u(\chi_u+C_0\omega)^{m+1}\wedge (\chi+C_0\omega)^l\wedge\omega^{n-m-l-1}\\
				&-\sum_{p=0}^m\sum_{q=0}^lC_{pq}\int_Mu\chi_u^p\wedge\chi^q\wedge \omega^{n-p-q}
			\end{aligned}
		\end{equation*}
		for some constants $C_{pq}$.
		Now
		\begin{align}\label{ReCY2}
			&\left|\int_M{u\ \left(\chi_u+C_0\omega\right)^{m+1}\wedge  (\chi+C_0\omega)^{l}\wedge\omega^{n-m-1-l} }\right|\notag\\
			\le &|u|_{L^{\infty}}\left|\int_M{ \left(\chi_u+C_0\omega\right)^{m+1}\wedge  (\chi+C_0\omega)^{l}\wedge\omega^{n-m-l-1} }\right|\notag\\
			=&|u|_{L^{\infty}}\left|\int_M{  (\chi+C_0\omega)^{m+l+1}\wedge\omega^{n-m-l-1} }\right|\notag\\
			\le& C_1|u|_{L^{\infty}}
		\end{align}
		and then by inductive assumption, inequality (\ref{ReCY1}) follows.
	\end{proof}
	
	\section{The existence of the longtime solution and proof of Theorem \ref{els}}
	In this section we prove Theorem \ref{els}, i.e. the following
	
	\begin{theorem}\label{3longtimeexistence}
		Let $(M,\omega)$ be a compact K\"{a}hler manifold and $\chi$  a closed real $(1,1 )$ form with $\theta_0\in(0, \pi)$.
		If $u_0$ satisfies inequality \eqref{B0}, then  flow \eqref{FYZflow} has  a  unique smooth longtime solution $u$.
	\end{theorem}
	
	We assume that $u$ is the solution of our flow \eqref{FYZflow} in $M\times[0, T)$, where $T$ is the maximal existence time. By showing the uniform a priori estimates in the following subsections, we can prove $T=\infty$.
	
	\subsection{The $u_t$ estimate }
	
	\begin{lemma}\label{utestimate}
		Let $u(x, t)$ be a solution of  flow \eqref{FYZflow} with the initial data satisfying  \eqref{B0}.
		For any $(x, t)\in M\times [0, T)$,
		\begin{align}\label{060303}
			\min_M u_t|_{t=0}\leq u_t(x,t)\leq \max_M u_t|_{t=0};
		\end{align}
		in particular,
		\begin{align}\label{060902}
			0<\min\limits_{M}\theta_{\omega}(\chi_{ u_0(x)})\le \theta_{\omega}(\chi_{u(x,t)})\le B_0(u_0)<\pi.
		\end{align}
	\end{lemma}
	\begin{proof}
		The  $u_t$ satisfies the equation:
		\begin{align*}
			(u_t)_t=F^{i\bar j} (u_t)_{i\bar j}.
		\end{align*}
		By the maximum principle, $u_t$ attains its maximum and minimum on the initial time, i.e.,
		inequality (\ref{060303}) holds, i.e.,
		\begin{align*}
			\min_M\cot\theta_{\omega}(\chi_{ u_0}) \le u_{t}(x,t)+\cot\theta_0\le\max_M \cot\theta_{\omega}(\chi_{u_0}),
		\end{align*}
		or
		\begin{align*}
			\min\limits_{ M}\cot\theta_{\omega}(\chi_{{u_0}})
			\leq  \cot \theta_{\omega}(\chi_{u(x,t)})\le \max\limits_{ M}\cot\theta_{\omega}(\chi_{{u_0}}).
		\end{align*}
		Thus we obtain
		\begin{align*}
			0<\min\limits_{M}\theta_{\omega}(\chi_{u_0})\le\theta_{\omega}(\chi_{u(x,t)})\le \max\limits_{M}\theta_{\omega}(\chi_{u_0})=B_0(u_0).
		\end{align*}
	\end{proof}
	We have a useful corollary of the above lemma.
	
	\begin{corollary}\label{060904}
		Let   $\lambda_n(x, t)$ be the smallest eigenvalue of $\chi_{u}$ with respect to the metric $\omega$ at $(x, t)$. Then
		\begin{align*}
			\max\limits_{M\times[0, T)}|\lambda_n|\le A_1\ \textup{for \
				$A_1:=|\cot B_0(u_0)|+\Bigl|\cot\Bigl(\frac{\min_{M}\theta_{\omega}(\chi_{u_0})}{n}\Bigr)\Bigr|$}.
		\end{align*}
	\end{corollary}
	
	\begin{proof}
		By Lemma \ref{utestimate}, we have
		\begin{align*}
			0<\frac{\min_M\theta_{\omega}(\chi_{u_0})}{n}\le \frac{\theta_{\omega}(\chi_u)}{n}\le\arccot\lambda_n\le B_0(u_0)<\pi.
		\end{align*}
		Hence  we have
		\begin{align*}
			\cot B_0(u_0)\le \lambda_n\le \cot\Bigl(\frac{\min_{M}\theta_{\omega}(\chi_{ u_0})}{n}\Bigr).
		\end{align*}
	\end{proof}
	\subsection{The complex Hessian estimate}
	For any $T_0<T$,
	we have proved $u_t$ is uniformly bounded and thus $|u|\le CT_0+|u_0|_{C^0}$ in $M\times[0, T_0]$. We  next prove the complex Hessian estimate.
	
	\begin{proposition}
		Let $u(x, t)$ be a solution of  flow \eqref{FYZflow} with the initial data satisfying  \eqref{B0}.
		There exists a uniform constant\ $C$ such that
		\begin{align*}
			\sup\limits_{M\times[0, T_0]}|\partial \bar\partial u|_{\omega}\le Ce^{CT_0}.
		\end{align*}
	\end{proposition}
	
	\begin{proof}
		Denote\ $w_{i\bar j}:=\chi_{i\bar j}+u_{i\bar j}$ as before.  Denote $S(T^{1,0}M):=\bigcup_{x\in M}\{\xi\in T^{1,0}_xM\,|\,|\xi|_{\omega}=1\}$. Consider on $S(T^{1,0}M)\times [0,T_0]$
		the auxiliary function
		\begin{align*}
			\tilde Q(x, t,\xi(x))=\log(w_{i\bar j}\xi^{i}\bar{\xi}^{j})-K_0t,
		\end{align*}
		where $K_0$ is a uniformly large constant to be chosen later.
		
		Suppose the function $\tilde Q$ attains its maximum at $(x_0, t_0)$ along the direction $\xi_0=\xi(x_0)$. We will prove that  $t_0=0$ and thus the estimate follows. If $t_0>0$,  we choose holomorphic  coordinates near \ $x_0$ such that
		\begin{equation}\label{3coordinate}
			\begin{aligned}
				& g_{i\bar j}(x_0)=\delta_{i\bar j}, \ \partial_kg_{i\bar j}(x_0)=0, \text{and}\\
				&w_{i\bar j}(x_0, t_0)=\lambda_{i}\delta_{i\bar j}\ \textup{with}\ \lambda_1\ge \lambda_2\ge \ldots\ge \lambda_n
			\end{aligned}
		\end{equation}
		which forces \ $\xi_0=\frac{\partial }{\partial z_1}$. We extend\ $\xi_0$ near\ $x_0$ as\ $ \tilde\xi_0(x)=(g_{1\bar 1})^{-\frac{1}{2}}\frac{\partial}{\partial z_1}$. Then the function\ $Q(x, t)=\tilde Q(x, t,\tilde\xi_0(x))$ on $M\times [0,T_0]$ attains its maximum at \ $(x_0, t_0)$.
		
		By the maximum principle,   we have at $(x_0, t_0)$
		\begin{align*}
			0\le& Q_t=\frac{u_{t1\bar 1}}{w_{1\bar 1}}-K_0,\notag\\
			0=&Q_i=\frac{w_{1\bar 1, i}}{w_{1\bar 1}},\\
			0\leq &-Q_{i\bar i}=-\frac{w_{1\bar 1, i\bar i}}{w_{1\bar 1}}+\frac{|w_{1\bar 1,i}|^2}{w_{1\bar 1}^2}=-\frac{w_{1\bar 1, i\bar i}}{w_{1\bar 1}}.\notag
		\end{align*}
		Hence we have
		\begin{align}\label{3PQ1}
			0\le Q_t-F^{i\bar i}Q_{i\bar i}
			=\lambda_1^{-1}(u_{t1\bar 1}-F^{i\bar i}w_{1\bar 1, i\bar i})  -K_0.
		\end{align}
		Since $d\chi=0$, by covariant  derivative formulae, we have
		\begin{align}
			w_{1\bar 1, i\bar i}=w_{i\bar i, 1\bar 1}+(\lambda_1-\lambda_i)R_{1\bar1 i\bar i}.\label{3L060651}
		\end{align}
		On the other hand, by \eqref{060601}, we have
		\begin{align}
			u_{t1\bar 1}- F^{i\bar i}w_{i\bar i, 1\bar 1}=&- F^{i\bar i}(1+\lambda_j^2)^{-1}(\lambda_i+\lambda_j)|w_{i\bar j, 1}|^2\notag\\
			&+2\cot\theta_{\omega}(\chi_u)\csc^2\theta_{\omega}(\chi_u)\frac{w_{i\bar i, \bar 1}}{1+\lambda_i^2}\frac{w_{j\bar j, \bar 1}}{1+\lambda_j^2}\notag\\
			=&-\sum\limits_{i\neq j}F^{i\bar i}(\lambda_i+\lambda_j)\frac{|w_{i\bar j, 1}|^2}{1+\lambda_j^2}-2F^{i\bar i}\lambda_i\frac{|w_{i\bar i, 1}|^2}{(1+\lambda_i^2)^{2}}\notag\\
			&+2\cot\theta_{\omega}(\chi_u)\csc^2\theta_{\omega}(\chi_u)\frac{w_{i\bar i, \bar 1}}{1+\lambda_i^2}\frac{w_{j\bar j, \bar 1}}{1+\lambda_j^2}. \label{3L060801}
		\end{align}
		However since  $\cot\theta(\lambda)$ is concave, by  (\ref{060702})
		\begin{align}
			u_{t1\bar 1}-F^{i\bar i}w_{i\bar i,1\bar 1}\leq -\sum\limits_{i\neq j}F^{i\bar i}(1+\lambda_j^2)^{-1}(\lambda_i+\lambda_j)|w_{i\bar j, 1}|^2\le 0,\label{3L6651}
		\end{align}
		since  $\lambda_i+\lambda_j>0$ for any $i\neq j$.
		
		Inserting \eqref{3L060651} and \eqref{3L6651} into  \eqref{3PQ1}, we have
		\begin{align}\label{PQ0}
			0\le Q_t-F^{i\bar i}Q_{i\bar i}\le 2|Rm|_{C^0}\sum\limits_{i=1}^n F^{i\bar i}-K_0.
		\end{align}
		Noting that $\sin\theta_{\omega}(\chi_u)\ge \min\left\{\sin{B_0(u_0)},\sin{\left(\min_M\theta_{\omega}(\chi_{u_0})\right)}\right\}$, for any $1\le i\le n$ we have \begin{align*}
			F^{i\bar i}=&\frac{1}{\sin^2\theta_{\omega}(\chi_u)(1+\lambda_i^2)}\\
			\le& \frac{1}{\min\left\{\sin^2{B_0(u_0)},\sin^2{\left(\min_M\theta_{\omega}(\chi_{u_0})\right)}\right\}}:=A_2.
		\end{align*}
		Inserting the above into \eqref{PQ0} and choosing $K_0=2nA_2|Rm|_{C^{0}}+1$, we have
		\begin{align}
			0\le Q_t-F^{i\bar i}Q_{i\bar i}\le 2nA_2|Rm|_{C^0}-K_0=-1,
		\end{align}
		which is a contradiction.
		Therefore $t_0=0$ and then for any $t\in[0,T_0]$, it holds
		
		\begin{align*}
			w_{i\bar j}\xi^{i}\bar{\xi}^{j}(x,t) e^{-K_0t}\le w_{1\bar 1}(x,0) =
			{u(0)}_{1\bar 1}+\chi_{1\bar 1}\le C.
		\end{align*}
	\end{proof}
	
	\subsection{\textbf{Proof of  Theorem \ref{3longtimeexistence}}}
	Since we have proved the $u_t$ estimate, the $C^0$ estimate and the complex Hessian estimate, by the concavity of the flow  \eqref{FYZflow}, we can apply the Evans-Krylov theory to get the higher order estimates of the solution.
	
	If the   maximal existence time $T<\infty$, then $u$ is uniformly $C^{k}$-bounded (for any $k\ge 0$) in $M\times[0, T]$ and then there exists $\epsilon>0$ such that the flow exists on $M\times[0, T+\epsilon_0]$ , which  is a contradiction since $T$ is the maximal existence time. Thus $T=\infty$.
	
	\section{Convergence of longtime solution and proof of Theorem \ref{mainthm}}
	In this section, we prove Theorem \ref{mainthm}, i.e., the following
	
	\begin{theorem}\label{mainthm1}
		Let $(M,\omega)$ be a compact K\"{a}hler manifold of dimension $n$ and $\chi$  a closed real $(1,1 )$ form with $\theta_0\in(0, \pi)$. Suppose the LYZ equation \eqref{LYZ2} has a subsolution  $\underline{u}$   in the sense of \eqref{A0} which also satisfies \eqref{B01}.
		If $u_0$ satisfies \eqref{B0}, then the longtime solution $u(x,t)$ of  flow \eqref{FYZflow} converges to a smooth solution $u^{\infty}$ to the  LYZ equation:
		\begin{align*}
			\theta_{\omega}(\chi_{u^{\infty}})=\theta_0.
		\end{align*}
	\end{theorem}
	
	\subsection{The $C^0$ estimate}
	We first prove a Harnack type inequality along our flow.
	\begin{lemma}\label{Harnack}
		Let $u$ be the solution of flow \eqref{FYZflow} on $M\times [0,\infty)$. Then for any $T_0<\infty$ we have the following Harnack type inequality:
		\begin{align*}
			\sup\limits_{M\times[0, T_0]}u(x,t)\le C\Bigl(-\inf\limits_{M\times[0, T_0]}\bigl(u(x,t)- u_0(x)\bigr)+1\Bigr).
		\end{align*}
	\end{lemma}
	
	\begin{proof}
		For any $t\in[0, T_0]$,  we have  $\theta_{\omega}(\chi_{u(t)})\leq B_0(u_0)<\pi$ by Lemma \ref{utestimate}.
		Then by the convexity of $\Gamma_{\omega,B_0(u_0)}:=\{\alpha\in\Lambda^{1,1}(M,\mathbb R)\ |\\
		\ \theta_{\omega}(\alpha)<B_0(u_0)\}$
		in Lemma \ref{convexity}, we have
		\begin{align*}
			\theta_{\omega}(\chi_{su+(1-s)u_0})\leq B_0(u_0)<\eta_0<\pi,
		\end{align*}
		where $\eta_0=B_0(u_0)/6+5\pi/6$ for convenience. Hence,
		\begin{align}\label{us}
			&\frac{\textup{Im}(\chi_{su(t)+(1-s)u_0}+\sqrt{-1}\omega)^n}{\omega^n}\notag\\
			=&\prod_{k=1}^n(1+\lambda_k^2(\chi_{su(t)+(1-s) u_0}))^{\frac 1 2}\sin\theta_{\omega}(\chi_{su(t)+(1-s) u_0})\notag\\
			\geq&\left\{
			\begin{aligned}
				&\sin\,\eta_0,\ \textup{if $\theta_{\omega}(\chi_{su(t)+(1-s)u_0})\geq \frac\pi 6$}\\
				&\sqrt{1+\lambda_1^2}\sin\arccot\lambda_1=1,\ \textup{if $\theta_{\omega}(\chi_{su(t)+(1-s)u_0})< \frac\pi 6$}
			\end{aligned}
			\right.\notag\\
			\geq &2c_0:=\sin{\eta_0}.
		\end{align}
		
		By Lemma \ref{IMCY},  the imaginary part of the Calabi-Yau functional is constant along the flow. Hence,
		\begin{align}\label{053101}
			0=&\mathrm{Im}\bigl(\textup{CY}_{\mathbb{C}}(u(t))\bigr)-\mathrm{Im}\bigl(\textup{CY}_{\mathbb{C}}(u_0)
			\bigr)\notag\\
			=&\int_{0}^{1}{\frac{d}{ds}\mathrm{Im}\bigl(\textup{CY}_{\mathbb{C}}(s u(t)+(1-s)u_0)\bigr)ds}\notag\\
			=&\int_{0}^{1}{\int_{M}{ (u(t)-u_0)}\mathrm{Im}\bigl(\chi_{s u(t)+(1-s)u_0}+\sqrt{-1}\omega\bigr)^n}ds\notag\\
			=&{\int_{M}(u(t)-u_0)\Bigl(\int_{0}^{1}\mathrm{Im}\bigl(\chi_{s u(t)+(1-s){u_0}}+\sqrt{-1}\omega\bigr)^nds\Bigr)}.
		\end{align}
		Thus we have
		\begin{align*}
			&\int_M { (u-{u_0})   \omega^n }\notag\\
			=&\int_M { (u-u_0)
				\omega^n }-\frac{1}{c_0}\int_{M}(u-u_0)\Bigl(\int_{0}^{1}\mathrm{Im}\bigl(\chi_{su(t)+(1-s)u_0}
			+\sqrt{-1}\omega\bigr)^nds\Bigr)\\
			=&\frac{1}{c_0}\int_M  -(u-u_0) \underbrace{\Bigl( -c_0 \omega^n +\int_{0}^{1}\mathrm{Im}\bigl(\chi_{su(t)+(1-s)u_0}+\sqrt{-1}\omega\bigr)^nds \Bigr)}_{\text{This term is nonnegative by} \ \eqref{us}  } \\
			\le&\frac{-\inf\limits_{M\times[0,T_0]}(u-u_0)}{c_0}\int_M\Bigl(-c_0 \omega^n  +\int_{0}^{1}\mathrm{Im}{\bigl(\chi_{su(t)+(1-s)u_0}
				+\sqrt{-1}\omega\bigr)^n}ds \Bigr)\\
			=&\frac{-\inf\limits_{M\times[0,T_0]}(u-u_0)}{c_0}\Bigl(-c_0\int_M \omega^n  +\int_{0}^{1}\mathrm{Im}\int_M{\bigl(\chi_{su(t)+(1-s)u_0}
				+\sqrt{-1}\omega\bigr)^n}ds \Bigr)\\
			=&\frac{-\inf\limits_{M\times[0,T_0]}{(u-u_0)}}{c_0}\Bigl(-c_0\int_M \omega^n  +\mathrm{Im}\int_M{\bigl(\chi
				+\sqrt{-1}\omega\bigr)^n} \Bigr)\\
			\le& c_0^{-1}\mathrm{Im}\int_M{\bigl(\chi
				+\sqrt{-1}\omega\bigr)^n} \Bigl(-\inf\limits_{M\times[0,T_0]}(u-u_0)\Bigr)\\
			=&C \Bigl(-\inf\limits_{M\times[0,T_0]}(u-u_0)\Bigr),
		\end{align*}
		where $C=c_0^{-1}\mathrm{Im}\int_M{(\chi
			+\sqrt{-1}\omega)^n}$. Therefore  we have
		\begin{align}\label{harnack2}
			\int_M u(x,t)\omega^n\le C\Bigr(-\inf_{M\times[0,T_0]}\bigl(u(x,t)-u_0(x)\bigr)+1\Bigr).
		\end{align}
		
		On the other hand, let $G(x, z)$ be Green's function of the metric $\omega$ on $M$.
		Then for any $(x, t)\in M\times[0, T_0]$,
		\begin{align*}
			u(x, t)=\Bigl(\int_M\omega^n\Bigr)^{-1}\int_{M} u(z, t)\omega^n-\int_{z\in M}\Delta_{\omega} u(z, t) G(x, z)\omega^n.
		\end{align*}
		Since $\Delta_{\omega} u>-\mathrm{tr}_{\omega}\chi>-C_0$ and $G(x,y)$ is bounded from below,  there exits a uniform constant $C$ such that
		\begin{align}\label{harnack3}
			u(x, t)\le \Bigl(\int_M\omega^n\Bigr)^{-1}\int_{M} u(z, t)\omega^n+C.
		\end{align}
		
		Combing \eqref{harnack2} with \eqref{harnack3}, we obtain the desired estimate.
	\end{proof}
	Now we can prove the $C^0$ estimate similar as Phong-T\^o \cite{pt2017arxiv}.
	\begin{proposition}\label{0606zoero}
		Along flow \eqref{FYZflow}, there exists a uniform constant $M_0$ independent of $T$ such that
		\begin{align*}
			|u|_{C^0(M\times[0,\infty))}\le M_0.
		\end{align*}
	\end{proposition}
	\begin{proof}
		Combining (\ref{053101}) with (\ref{us}) implies for any $t\in [0,\infty)$,
		\begin{align*}
			\sup\limits_{x\in M}(u(x,t)-u_0(x))\ge 0.
		\end{align*}
		Combing the above inequality with the concavity of the equation, we can apply  Lemma 1 by  Phong-T\^o \cite{pt2017arxiv}: there exists a uniform constant  $C_1$ such that
		\begin{align*}
			\inf_{M\times [0,T_0]}(u-\underline{u})\ge -C_1\quad \textup{ for any $T_0>0$}.
		\end{align*}
		Then combing this estimate with the  Harnack type inequality in Lemma \ref{Harnack}, we have
		\begin{align*}
			\sup\limits_{M\times[0, T_0]}u\le C.
		\end{align*}
		Since $T_0$ is arbitrary, the result follows.
	\end{proof}

	\subsection{The gradient estimate}
	We can use the following lemma by  Phong-T\^o which plays an important role in the gradient and second order estimates. In fact, it follows from the concavity of the function $\cot \theta(\chi_{{u}})$.
	
	\begin{lemma}\label{GJ}\cite{pt2017arxiv}
		Let $\delta$ and $K$ be two constants in Definition \ref{060901}.
		There exists a constant $\kappa_0$ depending only on  \ $\delta$, $K$, $\underline{u},\, (M, \omega)$,
		and $\chi$ such that if
		\begin{equation*}
			1+\lambda_1^2>\max\Bigl\{(K+\max\limits_{M}|\lambda(\chi_{\underline{u}})|+1)^2, \kappa_0^{-1}(1+A_1^2)\Bigr\},
		\end{equation*}
		then
		\begin{align}
			u_t-\sum F^{i\bar j}(u_{i\bar j}-\underline{u}_{i\bar j})\ge \kappa_0\sum F^{i\bar i}.
		\end{align}
	\end{lemma}
	
	We prove the gradient estimate following the argument in the elliptic case by Collins-Yau \cite{cy2018arxiv}.

	\begin{proposition}\label{0606gradient}
		Let $u$ be the solution of  flow \eqref{FYZflow}. There exists a uniform constant $M_1$ such that
		\begin{align*}
			\max\limits_{M\times[0,\infty)}|\nabla u|_{\omega}\le M_1.
		\end{align*}
	\end{proposition}
	
	\begin{proof}
		Without loss of generality, we assume $\underline u=0$; otherwise we write $\chi_u=\chi_{\underline u}+i\partial\bar\partial (u-\underline u)$  and replace $\chi$ by $\chi_{\underline u}$ and $u$ by $u-\underline u$.
		
		We consider the function
		\begin{align*}
			\tilde G=|\nabla u|^2\exp{\psi(u)}
		\end{align*}
		where
		\begin{align*}
			\psi(u)=-D_0 u+\bigl({ u+M_0+1}\bigr)^{-1}
		\end{align*}
		where $M_0$ is from Proposition \ref{0606zoero} and $D_0$ is a constant to be determined later.
		
		For any fixed time $T_0< \infty$, assume the function $\tilde G$ on $M\times [0,T_0]$ attains its maximum
		at $(x_0, t_0)$.  If\ $t_0=0$, we have the desired estimate directly. Hence  we assume\ $t_0>0$.
		The function $G:=\log \tilde G=\log|\nabla u|^2+\psi(u)$ also attains its maximum at $(x_0, t_0)$.
		By the maximum principle, we have  $\mathcal{P}G(x_0, t_0)\ge 0$ .
		
		Take the holomorphic coordinates \eqref{3coordinate} near $x_0$.
		By \eqref{Fijdef}
		\begin{align*}
			F^{i\bar j}(x_0, t_0)=\frac{\csc^2\theta(\lambda)}{1+\lambda_i^2}\delta_{ij}.
		\end{align*}
		
		We take  the manipulation at $(x_0, t_0)$:
		\begin{align*}
			G_t=&\frac{u_{kt}u_{\bar k}+u_{k}u_{\bar k t}}{|\nabla u|^2}+\psi' u_{t},\\
			G_{i}=&\frac{u_{ki}u_{\bar k}+u_{k}u_{\bar k i}}{|\nabla u|^2}+\psi' u_{i}=0,
			\\
			G_{i\bar j}=&\frac{u_{ki\bar j}u_{\bar k}+u_{ki}u_{\bar k\bar j}+u_{k\bar j}u_{\bar k i}+u_{k}u_{\bar k i\bar j}}{|\nabla u|^2}\\
			&-\frac{(u_{ki}u_{\bar k}+u_{k}u_{\bar k i})(u_{l\bar j}u_{\bar l}+u_{l}u_{\bar l \bar j})}{|\nabla u|^4}+\psi' u_{i\bar j}+\psi'' u_{i} u_{\bar j}.
		\end{align*}
		Hence
		\begin{align}\label{0Grad1}
			0\le& \mathcal{P}G=G_t-F^{i\bar i}G_{i\bar i}\notag\\
			=&{\frac{(u_{kt}-F^{i\bar i}u_{ki\bar i})u_{\bar k}+(u_{\bar k t}-F^{i\bar i}u_{\bar k i\bar i})u_k}{|\nabla u|^2}}\ (\textup{denoted by (I)})\\
			&-\frac{F^{i\bar i}(u_{ki}u_{\bar k\bar i}+u_{k\bar i}u_{\bar k i})|\nabla u|^2-F^{i\bar i}|\nabla_i|\nabla u|^2|^2}{|\nabla u|^4}\ (\textup{denoted by (II)})\notag\\
			&+\psi' (u_t-F^{i\bar i} u_{i\bar i})-\psi''F^{i\bar i} | u_{i}|^2\notag.
		\end{align}
		
		We first estimate   term (I).  By covariant derivatives formula and (\ref{060600}), we have
		\begin{align*}
			\textup{(I)}\leq&\frac{(u_{tk}-F^{i\bar i}u_{i\bar ik})u_{\bar k}+(u_{t\bar k}-F^{i\bar i}u_{ i\bar i\bar k})u_k+2F^{i\bar i}|Rm||\nabla u|^2}{|\nabla u|^2}\\
			\le &\frac{(u_{tk}-F^{i\bar i}w_{i\bar i,k}){u_{\bar k}}+(u_{t\bar k}-F^{i\bar i}w_{i\bar i, \bar k}){u_k}}{|\nabla u|^2}\\
			&+\frac{F^{i\bar i}(|\nabla\chi|+2|Rm||\nabla u|))}{|\nabla u|}\\
			=&\frac{F^{i\bar i}(|\nabla\chi|+2|Rm||\nabla u|))}{|\nabla u|}
			\le C_2.
		\end{align*}
		We then deal with  term (II). 
		Since\ $G_i=0$ for each $1\leq i\leq n$, we have
		\begin{align*}
			|\nabla_i|\nabla u|^2|^2=&{\bigl|\sum u_{ki}u_{\bar k}+\sum u_{k}u_{\bar k i}\bigr|^2}\\
			=&\bigl|\sum u_{ki}u_{\bar k}\bigr|^2+\bigl|\sum u_{k}u_{\bar k i}\bigr|^2+2\textup{Re}(\sum u_{ki}u_{\bar k}\sum u_{\bar k}u_{k\bar i  })\\
			=&\bigl|\sum u_{ki}u_{\bar k}\bigr|^2+\bigl|\sum u_{k}u_{\bar k i}\bigr|^2\\
			&+2\textup{Re}\Bigl(-\Bigl(\sum u_{k}u_{\bar k i}+|\nabla u|^2\psi'  u_{i}\Bigr)\sum u_{\bar k}u_{k\bar i  }\Bigr)\\
			=&\bigl|\sum u_{ki}u_{\bar k}\bigr|^2-\bigl|\sum u_{k}u_{\bar k i}\bigr|^2-2|\nabla u|^2\psi'\textup{Re} (u_{i}\sum u_{\bar k}u_{k\bar i }).
		\end{align*}
		Hence
		\begin{align*}
			\textup{(II)}
			=&-{|\nabla u|^{-2}}{F^{i\bar i}\Bigl(\sum |u_{ki}|^2+\sum |u_{k\bar i}|^2\Bigr)}+|\nabla u|^{-4}F^{i\bar i}\bigl|\sum u_{ki}u_{\bar k}\bigr|^2\\
			&-|\nabla u|^{-4}F^{i\bar i}\bigl|\sum u_{k}u_{\bar k i}\bigr|^2
			-2|\nabla u|^{-2}\psi'F^{i\bar i}\textup{Re} \Bigl( u_{i}\sum u_{\bar k}u_{k\bar i }\Bigr)\\
			\leq &-2|\nabla u|^{-2}\psi'F^{i\bar i}\textup{Re}(u_i\sum u_{\bar k}u_{k\bar i})
		\end{align*}
		where the last inequality holds by  the Cauchy-Schwarz inequality:
		\begin{align*}
			\bigl|\sum u_{ki}u_{\bar k}\bigr|^2\le \sum |u_{ki}|^2|\nabla u|^2.
		\end{align*}
		Since $u_{k\bar i}=w_{k\bar i}-\chi_{k\bar i}=\lambda_i\delta_{ki}-\chi_{k\bar i}$, by the Cauchy-Schwarz inequality again, we have
		\begin{align*}
			\textup{(II)}
			\leq &-2|\nabla u|^{-2}\psi'F^{i\bar i}|u_{i}|^2\lambda_{i }+2|\nabla u|^{-2}\psi'F^{i\bar i}\textup{Re} \Bigl( u_{i}\sum u_{\bar k} \chi_{k\bar i}\Bigr)\\
			\le&2|\psi'||\nabla  u|^{-1}\Bigl(\sum F^{i\bar i}| u_i|^2\Bigr)^{\frac{1}{2}}\Bigl(\sum F^{i\bar i}\lambda_i^2\Bigr)^{\frac{1}{2}}\\
			&+2|\chi||\psi'||\nabla u|^{-1}\Bigl(\sum F^{i\bar i}| u_i|^2\Bigr)^{\frac{1}{2}}\Bigl(\sum F^{i\bar i}\Bigr)^{\frac{1}{2}}.
		\end{align*}
		Clearly $\max\{\sum F^{i\bar i},\sum F^{i\bar i}\lambda_i^2\}\leq n\max_M\csc^2\theta_{\omega}(\chi_{u_0})$ by (\ref{060902}). \\If we take
		\begin{align}\label{06071}
			C_3:=4n\max_M\limits\csc\theta_{\omega}(\chi_{u_0})(1+\max\limits_M|\chi|),
		\end{align}
		then
		\begin{align*}
			\textup{(II)}\leq & C_3|\psi'||\nabla  u|^{-1}\Bigr(\sum F^{i\bar i}| u_i|^2\Bigr)^{\frac 1 2}.
		\end{align*}
		
		Inserting the estimates of (I) and (II) into \eqref{0Grad1}, we obtain
		\begin{align}\label{Grad2}
			0
			\le&\mathcal{P}G\le -\psi'(-u_t+F^{i\bar i}  u_{i\bar i})-\psi'' F^{i\bar i}| u_i|^2\notag\\
			&+C_3|\psi'||\nabla  u|^{-1}( F^{i\bar i}| u_i|^2)^{\frac 1 2}+C_2.
		\end{align}
		We  use the argument of Collins-Yau \cite{cy2018arxiv} and consider the two cases.
		Let $\epsilon_0$ be a positive constant satisfying
		\begin{align}\label{0607eps}
			\epsilon_0<\min\Bigl\{(K+\max\limits_{M}|\lambda(\chi_{\underline{u}})|+1)^{-1}, \kappa_0^{1/2}(1+A_1^2)^{-1/2},\notag\\ \frac{1}{2}C_3^{-1}\kappa_0(1+A_1^2)^{-1}\Bigr\}.
		\end{align}
		
		\textbf{Case 1: }\ $\sum\limits_{i=1}^nF^{i\bar i}| u_i|^2\ge \epsilon_0^2 |\nabla  u|^2$.
		
		By the definition of   $\psi$, $D_0\leq -\psi'\leq D_0+1$ and $\psi''=2(u-\inf_M u+1)^{-3}$. Hence,  by (\ref{Grad2})
		\begin{align*}
			0\le&- \frac{2\epsilon_0^2 |\nabla  u|^2}{( u+M_0 +1)^{3}}
			+(D_0+1)\Bigl(|u_t|+\frac{\csc^{2}\theta(\lambda)}{1+\lambda_i^2}|\lambda_i-
			\chi_{i\bar i}|\Bigr)\notag\\
			&+C_3(D_0+1)\csc\theta(\lambda)+C_2\\
			\le & - \frac{2\epsilon_0^2 |\nabla{ u}|^2}{( u+M_0 +1)^{3}}+ C(D_0+1).
		\end{align*}
		Thus we obtain
		\begin{align}\label{0607G1}
			|\nabla  u|^2\le C(D_0+1)\epsilon_0^{-2}( u+M_0+1)^{3}.
		\end{align}
		
		\textbf{Case 2:}\   $\sum\limits_{i=1}^nF^{i\bar i}| u_i|^2\le \epsilon_0^2 |\nabla  u|^2$.
		
		In this case, since $\psi''>0$, inequality (\ref{Grad2}) implies
		\begin{align}\label{060903}
			0\leq -\psi'(-u_t+F^{i\bar i}u_{i\bar i})+C_3(-\psi')\epsilon_0+C_2.
		\end{align}
		On the other hand, since $F^{1\bar 1}\le F^{i\bar i}$, we have
		\begin{equation*}
			\epsilon_0^2 | \nabla u|^2\ge  F^{1\bar 1}|\nabla  u|^2=\csc^2\theta(\lambda)\frac{|\nabla u|^2}{1+\lambda_1^2}.
		\end{equation*}
		Hence we get
		\begin{align*}
			{1+\lambda_1^2}\ge& \epsilon_0^{-2}\csc^2\theta(\lambda)\notag\\
			\ge& \epsilon_0^{-2}>\max\Bigl\{(K+\max\limits_{M}|\lambda(\chi_{\underline{u}})|+1)^2, \kappa_0^{-1}(1+A_1^2)\Bigr\}.
		\end{align*}
		Now we  apply the key Lemma \ref{GJ} to get
		\begin{align*}
			u_t-F^{i\bar j} u_{i\bar j}\ge \kappa_0\sum\limits_{i=1}^n F^{i\bar i}.
		\end{align*}
		Combined with \eqref{060903}, we get
		\begin{align}\label{06051}
			0\le \psi'\kappa_0\sum F^{i\bar i}+C_3(-\psi')\epsilon_0+C_2.
		\end{align}
		
		Since $\sum F^{i\bar i}> F^{n\bar n}=\frac{\csc^2\theta_{\omega}(\chi_u)}{1+\lambda_n^2}\geq (1+A_1^2)^{-1}$ by Corollary \ref{060904},
		and  $\epsilon_0<\frac{1}{2}C_3^{-1}\kappa_0(1+A_1^2)^{-1}$ by \eqref{0607eps},
		the sum of one half of the first term  and the second term  in \eqref{06051} is non-positive. Hence if we choose $D_0>2\kappa_0^{-1}C_2(1+A_1^2)$, we obtain the following contradiction.
		\begin{align*}
			0\le& \frac{1}{2}\psi'\kappa_0\sum F^{i\bar i}+C_2
			\le -\frac{D_0}{2}\kappa_0(1+A_1^2)^{-1}+C_2<0.
		\end{align*}
		
		Therefore if we choose $\epsilon_0$ satisfying \eqref{0607eps} and $D_0=2\kappa_0^{-1}C_2(1+A_1^2)+1$,   we really obtain the desired estimate \eqref{0607G1}.
	\end{proof}
	
	\subsection{Second order estimates}
	In the elliptic case, Collins-Jacob-Yau \cite{cjy2020cjm} used an auxiliary function containing the gradient term which modifies the one in Hou-Ma-Wu \cite{hmw2010mrl}. Here our auxiliary function does not contain the gradient term.
	\begin{proposition}\label{042601}
		There exists a uniform constant\ $M_2$ such that
		\begin{align*}
			\sup\limits_{M\times[0, \infty)}|\partial \bar\partial u|_{\omega}\le M_2.
		\end{align*}
	\end{proposition}
	\begin{proof}
		Without loss of generality, we assume that $\underline u=0$. Denote\ $w_{i\bar j}:=\chi_{i\bar j}+u_{i\bar j}$ as before. For any fixed\ $T_0<\infty$, we consider the auxiliary function
		on $S(T^{1,0}M)\times [0,T_0]$:
		\begin{align*}
			\tilde H(x, t,\xi(x))=\log(w_{i\bar j}\xi^{i}\bar{\xi}^{j})+\psi( u)
		\end{align*}
		where $\psi( u)=-D_1 u+{ u^2}/{2}$ with $D_1$ to be determined later.  Recall $M_0$ is the uniform bound of $|u|$ in Lemma \ref{0606zoero}.  Hence we have
		\begin{align}\label{0606psi}
			-D_1- M_0\le \psi'\le  -D_1+ M_0\quad \text{and} \quad
			\psi''=1.
		\end{align}
		
		Suppose the function $\tilde H$ attains its maximum at $(x_0, t_0)$ along the direction $\xi_0=\xi(x_0)$. If $t_0=0$, the estimate clearly holds. Hence we assume  $t_0>0$.  Take holomorphic  coordinates \eqref{3coordinate} near \ $x_0$ which forces \ $\xi_0=\frac{\partial }{\partial z_1}$. Extend\ $\xi_0$ near\ $x_0$ as\ $ \tilde\xi_0(x)=(g_{1\bar 1})^{-\frac{1}{2}}\frac{\partial}{\partial z_1}$. Then the function\ $H(x, t)=\tilde H(x, t,\tilde\xi_0(x))$ on $M\times [0,T_0]$ attains its maximum at \ $(x_0, t_0)$.
		
		By the maximum principle,   we have at $(x_0, t_0)$
		\begin{align}
			0\le& H_t=\frac{u_{t1\bar 1}}{w_{1\bar 1}}+\psi'u_t,\notag\\
			0=&H_i=\frac{w_{1\bar 1, i}}{w_{1\bar 1}}+\psi' u_i,\label{4060802}\\
			0\leq &-H_{i\bar i}=-\frac{w_{1\bar 1, i\bar i}}{w_{1\bar 1}}+\frac{|w_{1\bar 1,i}|^2}{w_{1\bar 1}^2}-\psi' u_{i\bar i}- |u_i|^2.\notag
		\end{align}
		Hence we have
		\begin{align}\label{PH1}
			0\le& H_t-F^{i\bar i}H_{i\bar i}\notag\\
			=&\lambda_1^{-1}(u_{t1\bar 1}-F^{i\bar i}w_{1\bar 1, i\bar i}) +\lambda_1^{-2}F^{i\bar i}|w_{1\bar 1, i}|^2 \quad \textup{(denoted by (I))}\notag\\
			&-F^{i\bar i}| u_i|^2+\psi'(u_t-F^{i\bar i} u_{i\bar i}).
		\end{align}
		
		We begin to deal with  term (I). By (\ref{3L060651}), we have
		\begin{align}
			u_{t1\bar 1}-F^{i\bar i}w_{1\bar 1, i\bar i}=u_{t1\bar 1}-F^{i\bar i}w_{i\bar i,1\bar 1}-F^{i\bar i}(\lambda_1-\lambda_i)R_{1\bar 1 i\bar i}.\label{060651}
		\end{align}
		On the other hand,  by \eqref{3L060801} and  by  (\ref{060702})   since  $\cot\theta(\lambda)$ is concave, we have
		\begin{align}
			u_{t1\bar 1}-F^{i\bar i}w_{i\bar i,1\bar 1}\leq -\sum\limits_{i\neq j}F^{i\bar i}(1+\lambda_j^2)^{-1}(\lambda_i+\lambda_j)|w_{i\bar j, 1}|^2.\label{6651}
		\end{align}
		Since  $\lambda_i+\lambda_j>0$ for any $i\neq j$, the above inequality implies
		\begin{align}
			u_{t1\bar 1}-F^{i\bar i}w_{i\bar i, 1\bar 1}
			\le& -\sum\limits_{i=2}^n\frac{\lambda_1+\lambda_i}{1+\lambda_1^2}F^{i\bar i}|w_{i\bar 1,1}|^2\notag\\
			=&-\sum\limits_{i=2}^n\frac{\lambda_1+\lambda_i}{1+\lambda_1^2}F^{i\bar i}|w_{1\bar 1,i}|^2.\label{060652}
		\end{align}
		Using \eqref{060651}, \eqref{060652} and \eqref{4060802},  we can estimate $\mathrm{(I)}$ as follows.
		\begin{align}
			\mathrm{(I)}\le& -\lambda_1^{-1}\sum\limits_{i=2}^n\frac{\lambda_1+\lambda_i}{1+\lambda_1^2}F^{i\bar i}|w_{1\bar 1,i}|^2+\lambda_1^{-2}\sum\limits_{i=1}^nF^{i\bar i}|w_{1\bar 1, i}|^2+C_4\notag\\
			=&\lambda_1^{-2}\sum\limits_{i=2}^nF^{i\bar i}|w_{1\bar 1, i}|^2\frac{1-\lambda_1\lambda_i}{1+\lambda_1^2}+\lambda_1^{-2}F^{1\bar 1}|w_{1\bar 1, 1}|^2+C_4\notag\\
			=&\psi'^2\sum\limits_{i=2}^nF^{i\bar i}
			| u_i|^2\frac{1-\lambda_1\lambda_i}{1+\lambda_1^2}+\psi'^2F^{1\bar 1}| u_1|^2+C_4.\notag
		\end{align}
		By Lemma \ref{060905}, we have $\lambda_1\geq\dots\geq\lambda_{n-1}\geq \cot(B_0(u_0)/2)$, and without loss of generality  we  assume $\lambda_1>1/\cot(B_0(u_0)/2)$. Hence for $2\leq i\leq  n-1$, $1-\lambda_1\lambda_i<0$.  For $i=n$, since $|\lambda_n|\leq A_1$, we have
		\begin{equation*}
			\frac{1-\lambda_1\lambda_n}{1+\lambda_1^2}\leq \frac{1+A_1}{\lambda_1}.
		\end{equation*}
		Hence
		\begin{align}
			\textup{(I)}\le&\psi'^2F^{1\bar 1}|\nabla  u|^2+ \psi'^2(1+A_1)\lambda_1^{-1}F^{n\bar n}| u_n|^2+C_4.\label{PFI}
		\end{align}
		
		Inserting  \eqref{PFI} into \eqref{PH1}, we have
		\begin{align}\label{1PH2}
			0 \le& \bigl(-1+(1+A_1)\psi'^2\lambda_1^{-1} \bigl)F^{n\bar n}
			| u_n|^2+\psi'^2F^{1\bar 1}|\nabla  u|^2\notag\\
			&+\psi'(u_t-F^{i\bar i} u_{i\bar i})+C_4\notag\\
			\le&\bigl(-1+(1+A_1)(D_1+ M_0)^2\lambda_1^{-1}\bigl)F^{n\bar n}
			| u_n|^2\notag\\
			&+(D_1+ M_0)^2M_1^2\csc^2\theta(\lambda)(1+\lambda_1^2)^{-1}
			+\psi'(u_t-F^{i\bar i}u_{i\bar i})+C_4.
		\end{align}
		The first term is negative if we assume
		\begin{align}\label{0606second1}
			\lambda_1>2(1+A_1)(D_1+ M_0)^2.
		\end{align}
		We further  assume
		\begin{align}\label{0606second2}
			1+\lambda_1^2>\max\Bigl\{(K+\max\limits_{M}|\chi_{\underline{u}}|)|+1)^2, \kappa_0^{-1}(1+A_1^2)\Bigr\}.
		\end{align}
		Then by Lemma \ref{GJ}, we  have
		\begin{equation*}
			u_t-F^{i\bar i} u_{i\bar i}\ge \kappa_0\sum\limits_{i=1}^n F^{i\bar i}\geq \kappa_0\frac{\csc\theta(\lambda)}{1+\lambda_n^2}\geq \kappa_0\frac{\csc\theta(\lambda)}{1+A_1^2}.
		\end{equation*}
		Hence if $D_1>M_0$, \eqref{1PH2} becomes
		\begin{align*}
			0\le&(D_1+ M_0)^2M_1^2\csc^2\theta(\lambda)(1+\lambda_1^2)^{-1}\notag\\
			&-(D_1- M_0)\csc^2\theta(\lambda){\kappa_0}(1+A_1^2)^{-1}+C_4
		\end{align*}
		or
		\begin{align*}
			\bigl((D_1- M_0){\kappa_0}(1+A_1^2)^{-1}
			-C_4\bigr)(1+\lambda_1^2)\le (D_1+ M_0)^2M_1^2.
		\end{align*}
		We choose
		\begin{equation*}
			D_1=(1+C_4)\kappa_0^{-1}(1+A_1^2)+ M_0.
		\end{equation*}
		Then we have
		\begin{align}\label{0606second3}
			\lambda_1\le( D_1+ M_0)M_1.
		\end{align}
		Combining \eqref{0606second1}, \eqref{0606second2} and  \eqref{0606second3}, we have
		$\lambda_1<C$ and then can obtain the desired $C^2$ estimate.
	\end{proof}
	
	\subsection{Proof of Theorem \ref{mainthm1}}
	The proof is the similar as the one in Phong-T\^o \cite{pt2017arxiv}. We sketch it for completeness. We have proved the uniform a priori estimates up to the second order. By the concavity of $\theta_{\omega}(\chi_u)$, we have the uniform $C^{2, \alpha}$ estimates and then the higher estimates hold.
	
	Since $u_t$ is uniformly bounded, there exists a constant $C$ such that $v(x, t):=u_t(x, t)+C>0$. Since $v$ satisfies $v_t=(u_t)_{t}=F^{i\bar j }(u_t)_{i\bar j}=F^{i\bar j}v_{i\bar j}$ and $F^{i\bar j}$ is uniformly elliptic, we can apply the  differential Harnack inequality (Cao \cite{ca1985im} and Gill  \cite{gi2011cag}) to get the following estimates
	\begin{align}\label{harnack}
		\max_{M}u_t(\cdot, t)-\min_{M}u_t(\cdot, t)=\max_{M}v(\cdot, t)-\min_{M}v(\cdot, t)\le Ce^{-C^{-1}t},
	\end{align}
	where $C$ is a uniform constant.
	
	By Lemma \ref{IMCY} and inequality (\ref{us})  we know that for any $t\in(0, \infty)$, there exists
	a point $x_0(t)$ such that $u_{t}(x_0(t), t)=0$.
	Therefore, for any $(x, t)\in M\times(0, \infty)$, by \eqref{harnack}, we have
	\begin{equation*}
		|u_t(x, t)|=|u_{t}(x, t)-u_{t}(x_0(t), t)|\le Ce^{-C^{-1}t},
	\end{equation*}
	and thus $u(x, t)$ converges exponentially to a function $u^{\infty}$. By the uniform  $C^k$ estimates  of $u(x, t)$  for all $ k\in \mathbb{N} $,  $u(x, t)$  converges to $u^{\infty}$ in $C^{\infty}$ and $u^{\infty}$ satisfies
	\begin{align*}
		\theta_{\omega}(\chi_{u^{\infty}}):=\sum\limits_{i=1}^n\arccot\lambda_i(\chi_{u^{\infty}})=\theta_0.
	\end{align*}
	
	\section{The convergence result on  K\"{a}hler surface with the semi-subsolution condition}
	
	In this section, we consider the compact K\"{a}hler surface case
	when $\chi$ satisfies the semi-subsolution condition i.e. $\chi-\cot \theta_0\omega\ge 0$.  We prove Theorem \ref{2dtheorem}, i.e.,
	
	\begin{theorem}\label{2dtheorem1}
		Let $(M, \omega)$ be a compact K\"{a}hler surface and $\chi$  a closed real $(1,1)$ form.
		Assume  $\theta_0\in (0,\pi)$ and $\chi\geq\cot\theta_0\omega$.
		Then there  exist a finite number of curves $E_i$ of negative self-intersection on $M$ such that
		the solution $u(x,t)$ of   flow \eqref{FYZflow} converges to a bounded function $u^{\infty}$ in $C^{\infty}_{loc}\left(M\setminus  \cup_i E_i\right)$ as $t$ tends to $\infty$
		with the following properties.
		\begin{enumerate}
			\item $\chi+\sqrt{-1}\partial\bar\partial u^{\infty}-\cot{B_1}\omega$ is a K\"ahler current which is smooth on $M{\setminus } \cup_i E_i$;
			\item $u^{\infty}$ satisfies the LYZ equation on $M{\setminus } \cup_i E_i$
			\begin{align}\label{2d}
				\mathrm{Re}(\chi_{u^{\infty}}+\sqrt{-1}\omega)^2=\cot \theta_0 \mathrm{Im}(\chi_{u^{\infty}}+\sqrt{-1}\omega)^2;
			\end{align}
			\item  $\chi_{u(x,t)}$ converges to $\chi_{u^{\infty}}$ and ${u^{\infty}}$  satisfies \eqref{2d} on $M$ in the sense of currents .
		\end{enumerate}
	\end{theorem}
	
	Here $u_0$ is a function in $\mathcal H_{B_1}$  for any $B_1\in (\theta_0, \pi)$.
	If $\theta_0\in(0, \frac{\pi}{2})$, we have $0\in \mathcal{H}_{B_1}$  for any $B_1\in (2\theta_0,  \pi )$.
	If $\theta_0\in [\frac{\pi}{2}, \pi)$, we first show that the semi-subsolution condition implies the non-empty of  $\mathcal{H}_{B_1}$ for any $B_1\in (\theta_0, \pi)$.
	
	\begin{lemma}\label{2dHnonempety}
		Let $(M, \omega)$ be a compact K\"{a}hler surface. Assume $\chi\ge\cot \theta_0\omega$ and $\theta_0\in [\frac{\pi}{2},\pi)$. Then for any $B_1\in (\theta_0, \pi)$, there exists a smooth function $u$ such that $u\in \mathcal{H}_{B_1}$.
	\end{lemma}
	
	\begin{proof}
		Let $\chi_{\epsilon}:=\chi-\epsilon\omega$ with $\epsilon>0$ sufficiently small.
		Define $\theta_0(\epsilon)$ as the principal argument of $\int_{M}(\chi_{\epsilon}+\sqrt{-1}\omega)^2$.  Then by definition,
		\begin{equation*}
			\cot\theta_0({\epsilon})=\frac{ \int_M \mathrm{Re} (\chi_{\epsilon}+\sqrt{-1}\omega)^2}{\int_M \mathrm{Im} (\chi_{\epsilon}+\sqrt{-1}\omega)^2 }.
		\end{equation*}
		Since $\theta_0\in(0, \pi)$, for any  $\epsilon>0$ sufficiently small  we have $\theta_0({\epsilon})\in(0, \pi)$ and thus $\mathrm{Im}\int_{M}(\chi_{\epsilon}+\sqrt{-1}\omega)^2=2\int_M \chi_{\epsilon}\wedge \omega  >0$. By direct manipulation, we have
		\begin{align*}
			\cot\theta_0({\epsilon})
			=&\frac{ \int_M  (\chi^2-\omega^2+ \epsilon^2\omega^2 -2   \epsilon \chi\wedge \omega)}
			{2\int_M  (\chi\wedge \omega   -\epsilon\omega^2  )          }\\
			=&\cot\theta_0-\epsilon+\epsilon\left(\cot\theta_0-\frac{\epsilon}{2}\right)\frac{ \int_M  \omega^2}
			{\int_M  (\chi\wedge \omega   -\epsilon\omega^2  )  }\\
			<& \cot\theta_0-\epsilon.
		\end{align*}
		This shows $\chi_{\epsilon}\ge \cot\theta_0\omega-\epsilon\omega>\cot\theta_0({\epsilon}) \omega$. Hence by Jacob-Yau \cite{jy2017ma}  there exists a smooth function $u_{\epsilon}$ solving
		\begin{align*}
			\mathrm{Re} (\chi_{\epsilon}+\sqrt{-1}\p\bar \p {u_{\epsilon}}+\sqrt{-1}\omega)^2  =  \cot\theta_0({\epsilon}) \mathrm{Im} (\chi_{\epsilon}+\sqrt{-1}\p\bar \p {u_{\epsilon}}+\sqrt{-1}\omega)^2.
		\end{align*}
		Thus for any $B_1\in (\theta_0, \pi)$, we have
		\begin{align*}
			\theta_{\omega}(\chi_{u_\epsilon})<\theta_{\omega}(\chi_{\epsilon}+\sqrt{-1}\p\bar \p u_\epsilon)=\theta_0(\epsilon)<B_1,
		\end{align*}
		where $\epsilon$ is sufficiently small since $\theta_0(\epsilon)$ attends to $\theta_0$ as $\epsilon$ goes to $0$ .
	\end{proof}
	
	We will use the  following proposition  proved by Song-Weinkove \cite{sw2008cpam}.
	
	\begin{proposition}[Song-Weinkove \cite{sw2008cpam}]\label{2dsongweinkove}
		Let $M$ be a  K\"{a}hler surface with a K\"{a}hler class $\beta\in H^{1,1}(M, \mathbb{R})$. If $\alpha\in H^{1,1}(M, \mathbb{R})$
		satisfies $\alpha^2>0$ and $\alpha\cdot \beta>0$, then either $\alpha$ is K\"{a}hler or there exists a positive integer m, curves of negative self-intersection $E_i, 1\le i\le m$ and positive numbers $a_i, 1\le i\le m$ such that
		\begin{align*}
			\alpha-\sum\limits_{i=1}^m a_i[E_i]
		\end{align*}
		is a K\"{a}hler class.
	\end{proposition}
	
	Since $2\cot\theta_0[\chi]\cdot[\omega]={[\chi]^2-[\omega]^2}$, if we let $\tilde \chi=\chi-\cot\theta_0\omega$, then we have
	\begin{align}\label{2d1}
		[\tilde\chi]^2=[\chi]^2 -2\cot\theta_0[\chi]\cdot[\omega]+\cot^2\theta_0[\omega]^2 =(1+\cot^2\theta_0 )[\omega]^2>0.
	\end{align}
	Since $\tilde \chi\ge 0$, we also have
	\begin{align*}
		[\tilde \chi]\cdot[\omega]>0,
	\end{align*}
	otherwise $\tilde \chi\equiv 0$ which contradicts with \eqref{2d1}.
	Hence we can apply Proposition \ref{2dsongweinkove} to get that there exists a finite number $m\ge 0$, curves of negative self-intersection $E_i, 1\le i\le m$ and positive numbers $a_i, 1\le i\le m$ such that
	$[\tilde \chi]-\sum_{i=1}^m a_i[E_i]$ is a K\"{a}hler class.
	
	Let $h_i$ be the hermitian metric on  $ [E_i]$  and $s_i$ be a holomorphic section of  $[E_i]$ which  vanishes along $E_i$ to order 1.
	Define
	\begin{equation*}
		S:=\sum\limits_{i=1}^m a_i\log|s_i|_{h_i}^2 ,
	\end{equation*}
	then
	\begin{equation}\label{043001}
		\tilde \chi+\sqrt{-1}\p\bar \p S>0.
	\end{equation}
	
	Similar as the argument in Section 2 in \cite{flsw2014apde}
	which is based on \cite{egz2009jams}, \cite{ts1988ma} and \cite{zh2006imrn}, we get the following result.
	\begin{lemma}\label{2dmongeampere1}
		Let $(M, \omega)$ be a compact K\"{a}hler surface.
		Assume $\tilde \chi:=\chi-\cot \theta_0\omega\ge 0$ and $\theta_0\in (0,\pi)$.
		Then there exists a unique (up to adding a constant ) bounded $\tilde \chi$-PSH function $v$ on $M$ and $v\in C^{\infty}_{loc}\left(M{\setminus }\cup_i E_i\right)$
		satisfying
		\begin{align}\label{2dmongeampere2}
			(\tilde \chi+\sqrt{-1}\partial\bar \partial v)^2=\csc^2\theta_0 \omega^2,
		\end{align}
		in the sense of currents.
	\end{lemma}
	
	\subsection{The uniform $C^0$-estimate }
	We have proved the uniform $u_t$ estimate and thus along the flow we have \begin{align*}\theta_{\omega}(\chi_u)\in (\min_M \theta_{\omega}(\chi_{u(0)}),  B_1).
	\end{align*}
	\begin{proposition}\label{2dzeroorder}
		Assume the same conditions hold as in Theorem \ref{2dtheorem}. Then there exists a uniform constant $M_0$
		such that for any $(x, t)\in M\times [0, \infty)$
		\begin{align}
			|u(x, t)|\le M_0.
		\end{align}
	\end{proposition}
	
	\begin{proof}
		For any $T_0$, we will prove $\sup\limits_{M\times[0, T_0]}|u(x, t)|\le M_0$.
		We use similar auxiliary functions by Fang-Lai-Song-Weinkove \cite{flsw2014apde} for the J-flow and Takahashi  \cite{ta2021cvpde} for the LBMCF.
		
		We first prove the upper bound of $u$ using the solution $v$ in Lemma \ref{2dmongeampere1}. Consider
		\begin{align*}
			w_{\varepsilon}(x, t):=u-(1+\varepsilon)v+\varepsilon S-C_0\varepsilon t,
		\end{align*}
		where $C_0$ is a large constant to be determined later.
		Since $w_{\varepsilon}(x, t)$ is  upper semi-continuous  on $M\times[0, T_0]$ with $w_{\varepsilon}=-\infty$ in $\cup_i E_i$, $w_{\varepsilon}$ attains
		its maximum on $M\times[0, T_0]$ at $(x_0, t_0)$ with $x_0\in M{\setminus }\cup_i E_i$. Our goal is to show $t_0=0$.
		
		At $(x_0, t_0)$, we have
		\begin{align}
			0\ge \sqrt{-1}\partial \bar\partial w_{\varepsilon}=&\sqrt{-1}\partial \bar\partial u-(1+\varepsilon)\sqrt{-1}\partial \bar\partial v+\varepsilon\sqrt{-1}\partial \bar\partial S\notag\\
			=&\tilde \chi_u-(1+\varepsilon)\tilde \chi_v+\varepsilon (\tilde \chi+\sqrt{-1}\p\bar \p S) \notag\\
			\ge& \tilde \chi_u-(1+\varepsilon)\tilde \chi_v,\label{2duniform1}
		\end{align}
		where in the last inequality we use inequality (\ref{043001}).
		Let $\lambda=(\lambda_1, \lambda_2)$ and $\mu=(\mu_1, \mu_2)$ be the eigenvalues of $\chi_u(x_0, t_0)$ and  $\tilde \chi_u(x_0, t_0)$ with respect to the metric $\omega$ respectively. Then $\lambda_i=\mu_i+\cot \theta_0$. Without loss of generality, we assume $\lambda_1\ge \lambda_2$. By direct manipulation, we have
		\begin{align}	\frac{dw_{\varepsilon}}{dt}(x_0, t_0)=&\frac{du}{dt}(x_0, t_0)-C_0 \varepsilon=\cot\theta(\chi_{u}(x_0,t_0))-\cot\theta_0-C_0 \varepsilon\label{2duniform22}\\
			=&\frac{\lambda_1\lambda_2-1}{\lambda_1+\lambda_2}-\cot\theta_0-C_0 \varepsilon\notag\\
			=&\frac{\mu_1\mu_2-\csc^2\theta_0}{\lambda_1+\lambda_2}-C_0 \varepsilon \label{2duniform2}.
		\end{align}
		
		\textbf{Case  1:}  $\mu_1\ge 0$ and $\mu_2\ge 0$. By \eqref{2duniform1} and \eqref{2dmongeampere2}, we have
		\begin{align}
			\mu_1\mu_2\le& (1+\varepsilon)^2\frac{\tilde \chi_v^2}{\omega^2}= (1+\varepsilon)^2\csc^2\theta_0.\label{2duniform3}
		\end{align}
		Inserting \eqref{2duniform3} into  \eqref{2duniform2}, we obtain
		\begin{align}
			\frac{dw_{\varepsilon}}{dt}(x_0, t_0)\le
			&\frac{\csc^2\theta_0}{\lambda_1+\lambda_2}(2\varepsilon+\varepsilon^2)-C_0 \varepsilon\notag\\
			\le& \frac{3\csc^2\theta_0}{\cot{\frac{B_1}{2}}- \cot{{B_1}}}\varepsilon-C_0\varepsilon\notag\\
			=&-\varepsilon<0, \label{2duniform4}
		\end{align}
		where we use $\lambda_1+\lambda_2\ge \cot{\frac{B_1}{2}}- \cot{{B_1}}>0$ and choose $C_0=\frac{3\csc^2\theta_0}{\cot{\frac{B_1}{2}}- \cot{{B_1}}}+1$.
		
		\textbf{Case 2:}  $\mu_1\ge 0$ and $\mu_2\leq 0$.
		By \eqref{2duniform2}, $\frac{dw_{\varepsilon}}{dt}(x_0, t_0)<-C_0\varepsilon<0$.
		
		\textbf{Case 3:}  $\mu_1\le 0$ and $\mu_2\le 0$.
		Then $\lambda_1=\mu_1+\cot\theta_0\le\cot\theta_0 $ and we get
		$\cot\theta(\chi_u(x_0, t_0))=\lambda_1-\frac{1+\lambda_1^2}{\lambda_1+\lambda_2}
		<\cot\theta_0$. Thus by \eqref{2duniform22}, we have $\frac{dw_{\varepsilon}}{dt}(x_0, t_0)=\frac{du}{dt}(x_0, t_0)-C_0 \varepsilon
		<0$.

		
		From the above three cases, we conclude $\frac{dw_{\varepsilon}}{dt}(x_0, t_0)<0$ and thus $t_0=0$.
		Thus   for any $\varepsilon>0$ and $(x,t)\in (M\setminus\cup_i E_i)\times[0, T_0]$,we have
		\begin{equation*}
			u(x, t)\le {u_0}(x_0)-(1+\varepsilon)v(x_0)+\varepsilon S(x_0) +(1+\varepsilon)v(x)-\varepsilon S(x)+C_0\varepsilon t.
		\end{equation*}
		Fix $(x,t)\in (M\setminus\cup_iE_i)\times [0,T_0]$ and let $\varepsilon $ tend to $0$, since $S$ is upper bounded, we have $u(x,t) \le\max{u_0}+2\max  |v| $,
		which also holds for any $(x,t)\in M\times [0,T_0]$ by continuity of $u(x,t)$. Since $T_0$ is arbitrary, $u\le \max{u_0}+2\max  |v|$ in $M\times [0,\infty]$.
		
		Next we prove the lower bound of $u$.
		Consider
		\begin{align*}
			\tilde w_{\varepsilon}:=u-(1-\varepsilon)v-\varepsilon S+C_0\varepsilon t,
		\end{align*}
		where $C_0$ is a constant as above.
		Since $\tilde w_{\varepsilon}(x, t)$ is  lower semi-continuous  with $\tilde w_{\varepsilon}=+\infty$ in $\cup_i E_i$,  $\tilde w_{\varepsilon}$ attains
		its minimum in $M\times[0, T_0]$ at $(x_1, t_1)$ with $x_1\in M{\setminus } \cup_i E_i$.
		
		At $(x_1, t_1)$, we have
		\begin{align}
			0\le \sqrt{-1}\partial \bar\partial \tilde w_{\varepsilon}=&\sqrt{-1}\partial \bar\partial u-(1-\varepsilon)\sqrt{-1}\partial \bar\partial v-\varepsilon\sqrt{-1}\partial \bar\partial S\notag\\
			=&\tilde \chi_u-(1-\varepsilon)\tilde \chi_v-\varepsilon(\tilde \chi+\sqrt{-1}\partial \bar\partial S)\notag\\
			\le& \tilde \chi_u-(1-\varepsilon)\tilde \chi_v.\label{2duniform2.1}
		\end{align}
		This implies
		\begin{align*}
			\mu_1\mu_2\ge (1-\varepsilon)^2\frac{\tilde \chi_v^2}{\omega^2}=(1-\varepsilon)^2\csc^2\theta_0.
		\end{align*}
		Hence
		\begin{align*}
			\frac{d\tilde w_{\varepsilon}}{dt}(x_1, t_1)
			=&\frac{\mu_1\mu_2-\csc^2\theta_0}{\lambda_1+\lambda_2}+C_0 \varepsilon\\
			\ge&  -\frac{2\csc^2\theta_0}{\lambda_1+\lambda_2}\varepsilon+C_0\varepsilon>0.
		\end{align*}
		Thus $\tilde w_{\varepsilon}$ attains its minimum at $t_1=0$ and the lower bound of $u$ follows.
	\end{proof}
	
	Combining the above uniform estimate and Proposition \ref{recy} yields
	
	\begin{corollary}\label{2duniformboundReCY}Along the  flow,
		there exists a uniform constant  $C$ such that
		\begin{align}
			\mathrm{Re}{(\mathrm{CY}_{\mathbb{C}}(u))}\le C.
		\end{align}
	\end{corollary}
	
	\subsection{$C^k$-estimate in compact set $K\subset M\setminus\cup_i E_i$.}
	
	Since $\chi-\cot\theta_0\omega_0$ is only nonnegative, we could not apply Lemma \ref{GJ} directly. But we can prove a similar type inequality  as in Lemma \ref{GJ}. In fact,
	we consider $\tilde u=u-S$. Since $\chi- \cot \theta_0\omega\ge 0$ and all $E_i, 1\le i\le m$
	are negative self-intersection,
	we have $\chi-\cot\theta_0\omega+\sqrt{-1}\partial\bar\partial S>0$, and thus there exists a small constant $\delta>0$ such that
	\begin{align}\label{subsolution9.1}
		\chi+\sqrt{-1}\partial \bar \partial S\ge( \cot\theta_0+\delta )\omega.
	\end{align}
	We can prove the following useful inequality which is the key for us to prove the gradient estimate and the complex Hessian estimate.
	
	\begin{lemma}\label{2dsubsolutioninequality}
		Assume the same conditions hold as in Theorem  \ref{2dtheorem1}.
		There exist uniform constants $K_0>0$ and $c_0>0$ such that if $|\lambda(\chi_u)|>K_0$, then
		\begin{align*}
			u_t-F^{i\bar j}(u_{i\bar j}-S_{i\bar j})\ge c_0.
		\end{align*}
	\end{lemma}
	
	\begin{proof}
		Choose the normal coordinates at {{$(x, t)$}} as before. By \eqref{subsolution9.1} we have
		\begin{align}
			u_t-F^{i\bar j}(u_{i\bar j}-S_{i\bar j})
			=&\cot\theta_{\omega}(\chi_u)-\cot\theta_0-F^{i\bar i}(w_{i\bar i}-\chi_{i\bar i}-S_{i\bar i})\notag\\
			\ge &\cot\theta_{\omega}(\chi_u)-\cot\theta_0-F^{i\bar i}w_{i\bar i}\notag\\
			&+(\delta+\cot\theta_0)\sum_{i=1}^2F^{i\bar i}\label{subsubsolution}.
		\end{align}
		By \eqref{060904}, we have $|\lambda_2|\le A_1$. Recall that
		$\cot\theta_{\omega}(\chi_u)=\frac{\lambda_1\lambda_2-1}{\lambda_1+\lambda_2}$ and $\csc^2\theta_{\omega}(\chi_u)=1+\cot^2\theta_{\omega}(\chi_u)=\frac{(1+\lambda_1^2)(1+\lambda_2^2)}{(\lambda_1+\lambda_2)^2}$. Hence we have
		\begin{align}
			\cot\theta_{\omega}(\chi_u)-F^{i\bar i}w_{i\bar i}=&\frac{\lambda_1\lambda_2-1}{\lambda_1+\lambda_2}-\frac{(1+\lambda_1^2)\lambda_2}{(\lambda_1+\lambda_2)^2}-\frac{(1+\lambda_2^2)\lambda_1}{(\lambda_1+\lambda_2)^2}\notag\\
			=&\frac{-2}{\lambda_1+\lambda_2}
			\ge-C\lambda_1^{-1}\label{subsubsolution1}.
		\end{align}
		
		For the other terms in \eqref{subsubsolution}, we have
		\begin{align}
			-\cot\theta_0+& (\delta+\cot\theta_0)\sum_{i=1}^2F^{i\bar i}\\
			\ge& \cot\theta_0(\frac{\csc^2\theta_{\omega}(\chi_u)}{1+\lambda_2^2}-1)
			+\delta\sum_{i=1}^2F^{i\bar i}-C\lambda_1^{-1}\notag\\
			=&\cot\theta_0(\frac{1+\lambda_1^2}{(\lambda_1+\lambda_2)^2}-1)
			+\delta\sum_{i=1}^2F^{i\bar i}-C\lambda_1^{-1}\notag\\
			\ge& -C\lambda_1^{-1}+\delta\sum_{i=1}^2F^{i\bar i}\notag\\
			\ge& -C\lambda_1^{-1}+\delta\frac{(1+\lambda_1^2)}{(\lambda_1+\lambda_2)^2}\notag\\
			\ge& \frac{\delta}{2},\label{subsubsolution2}
		\end{align}
		where we assume $\lambda_1\ge K_0$ and choose $K_0$ sufficiently large.
		
		Inserting \eqref{subsubsolution1} and \eqref{subsubsolution2} into  \eqref{subsubsolution}, we obtain
		\begin{align*}
			u_t-F^{i\bar j}(u_{i\bar j}-S_{i\bar j})\ge\frac{\delta}{2}-C\lambda_1^{-1}
			\ge \frac{\delta}{4}.
		\end{align*}
	\end{proof}
	The following lemma is useful for us to prove the gradient estimate and the complex Hessian estimate.
	\begin{lemma}\label{2dSestimate}
		There exist uniform positive constants $\Lambda_0:=\min_{i}\{a_i^{-1}\}$ and $\Lambda_1$ such that for any $x\in M\setminus\cup_i E_i$, it holds
		\begin{align}\label{5.17}
			e^{\Lambda_0S(x)}\left(	|S(x)|^3+|\nabla S|^2(x)\right)\le \Lambda_1.
		\end{align}
	\end{lemma}
	\begin{proof}
		Since $S=\sum\limits_{i=1}^m a_i\log|s_i|_{h_i}^2$, there exists a uniform constant $C>0$  such that
		\begin{align}
			|\nabla S|^2\le C(1+\sum\limits_{i=1}^m|s_i|^{-2}).
		\end{align}
		Then we have \eqref{5.17} .
	\end{proof}
	\begin{proposition}\label{2dgradient1}
		There exist  uniform constants $D_0$ and $M_1$ such that for any $(x,t)\in {M\setminus\cup_i E_i\times[0, \infty)}$
		\begin{align}
			|\n
			u|_{\omega}(x,t)\le M_1{\prod_{i}|s_i|^{-D_0a_i}_{h_i}(x)}.
		\end{align}
	\end{proposition}
	
	\begin{proof}
		Since $S$ is upper semi-continuous, there exists a uniform constant $S_0$ such that $\sup_{M} S\le S_0$.
		We consider the function
		\begin{align*}
			G=\log|\nabla u|^2+{\psi( \tilde u)},
		\end{align*}
		where  $\tilde u=u-S$ and
		\begin{align*}
			\psi(\tilde u)=-D_0 \tilde u+\bigl({ \tilde u+S_0+M_0+1}\bigr)^{-1},
		\end{align*}
		where $D_0> \Lambda_0:=\min_{i}\{a_i^{-1}\}$ is a uniform constant to be determined later.
		
		Since $S$ is upper semi-continuous, we know that  $G$ is also upper semi-continuous.  Suppose that $\max\limits_{M\times [0, T_0]}G(x, t)=G(x_0, t_0)$. Since $S=-\infty$ on $\cup_i{E_i}$, we have $G(x, t)=-\infty$  on $\cup_i{E_i}$ and then $x_0\in M\setminus \cup_{i}E_i $.
		
		{If} $t_0=0$, we have for any $(x, t)\in M\setminus\cup_i E_i\times[0, T_0]$
		\begin{align}
			e^{G(x, t)}\le e^{G(x_0, 0)}\le |\nabla u_0|^2e^{D_0S_0+D_0M_0+S_0+M_0+1}\le M_{1,0},
		\end{align}
		where we used $S\le S_0$ and $M_{1,0}:=\max_M|\nabla u_0|^2e^{(D_0+1)(S_0+M_0)+1}$. This gives the estimate \eqref{2dgradient1}.
		
		In the following, we always assume  $t_0>0$.
		
		If $|\nabla u|(x_0, t_0)\le 2|\nabla S|(x_0, t_0)$, by Lemma \ref{2dSestimate}, we  get the desired estimate as follows
		\begin{align}
			e^{G(x_0, t_0)}\le& C|\nabla u|^2(x_0, t_0) e^{D_0 S(x_0)}\notag\\
			\le& 4C|\nabla S|^2(x_0, t_0) e^{D_0 S(x_0, t_0)}\le M_{1,1}.
		\end{align}
		Thus
		in the following, we always assume $|\nabla u|(x_0, t_0)\ge 2|\nabla S|(x_0, t_0)$ and then we have
		\begin{align}\label{2dgradientassumption}
			\frac{1}{2}	|\nabla u|(x_0, t_0)\le |\nabla\tilde u|(x_0, t_0)\le 2|\nabla u|(x_0, t_0).
		\end{align}
		Taking  the manipulation at $(x_0, t_0)$, we have
		\begin{align*}
			G_t=&\frac{u_{kt}u_{\bar k}+u_{k}u_{\bar k t}}{|\nabla u|^2}+\psi' u_{t},\\
			G_{i}=&\frac{u_{ki}u_{\bar k}+u_{k}u_{\bar k i}}{|\nabla u|^2}+\psi' \tilde{u}_{i}=0,
		\end{align*}
		and
		\begin{align}\label{Grad1}
			0\le& \mathcal{P}G=G_t-F^{i\bar i}G_{i\bar i}\notag\\
			=&{\frac{(u_{kt}-F^{i\bar i}u_{ki\bar i})u_{\bar k}+(u_{\bar k t}-F^{i\bar i}u_{\bar k i\bar i})u_k}{|\nabla u|^2}}\ (\textup{denoted by (I)})\notag\\
			&-\frac{F^{i\bar i}(u_{ki}u_{\bar k\bar i}+u_{k\bar i}u_{\bar k i})|\nabla u|^2-F^{i\bar i}|\nabla_i|\nabla u|^2|^2}{|\nabla u|^4}\ (\textup{denoted by (II)})\notag\\
			&+\psi' (u_t-F^{i\bar i} \tilde u_{i\bar i})-\psi''F^{i\bar i} | \tilde u_{i}|^2.
		\end{align}
		By the same estimate as that  in Proposition \ref{0606gradient}, we have
		\begin{align*}
			\textup{(I)}\leq C.
		\end{align*}
		We then deal with  term (II). 
		Since\ $G_i=0$ for each $1\leq i\leq 2$, we have
		\begin{align*}
			|\nabla_i|\nabla u|^2|^2
			=&\bigl|\sum u_{ki}u_{\bar k}\bigr|^2+\bigl|\sum u_{k}u_{\bar k i}\bigr|^2+2\textup{Re}\Bigl(\sum u_{ki}u_{\bar k}\sum u_{\bar k}u_{k\bar i  }\Bigr)\\
			=&\bigl|\sum u_{ki}u_{\bar k}\bigr|^2+\bigl|\sum u_{k}u_{\bar k i}\bigr|^2\notag\\
			&+2\textup{Re}\Bigl(-\Bigl(\sum u_{k}u_{\bar k i}+|\nabla u|^2\psi'  \tilde u_{i}\Bigr)\sum u_{\bar k}u_{k\bar i  }\Bigr)\\
			=&\bigl|\sum u_{ki}u_{\bar k}\bigr|^2-\bigl|\sum u_{k}u_{\bar k i}\bigr|^2-2|\nabla u|^2\psi'\textup{Re} \Bigl(\tilde u_{i}\sum u_{\bar k}u_{k\bar i }\Bigr).
		\end{align*}
		Hence
		\begin{align*}
			\textup{(II)}
			\leq &-2|\nabla u|^{-2}\psi'F^{i\bar i}\textup{Re}(\tilde u_i\sum u_{\bar k}u_{k\bar i}),
		\end{align*}
		Similar as the estimate in Proposition \ref{0606gradient}, we have
		\begin{align*}
			\textup{(II)}\leq & C|\psi'||\nabla  u|^{-1}\Bigr(\sum F^{i\bar i}| \tilde u_i|^2\Bigr)^{\frac 1 2}.
		\end{align*}
		Inserting the estimates of (I) and (II) into \eqref{Grad1}, we obtain
		\begin{align}\label{2dGrad2}
			0	\le\mathcal{P}G\le& -\psi'(-u_t+F^{i\bar i}  u_{i\bar i})-\psi'' F^{i\bar i}| \tilde u_i|^2\notag\\
			&+C|\psi'||\nabla  u|^{-1}( F^{i\bar i}|\tilde u_i|^2)^{\frac 1 2}+C.
		\end{align}
	We divide two cases to do the estimate.\\
	Let $\epsilon_0=\min\{\frac{1}{2}K_0^{-\frac{1}{2}}, \frac{c_0}{2C\min_{M}|\sin\theta_{\omega}(\chi_{u_0})|}\}$ where $K_0$ is the uniform constant in  Lemma \ref{2dsubsolutioninequality} and $C$ is the constant in \eqref{2dGrad2}.
	
	\textbf{Case 1: }\ $\sum\limits_{i=1}^2 F^{i\bar i}| \tilde u_i|^2\ge \epsilon_0^2 |\nabla  u|^2$.
	
	Since $D_0\leq -\psi'\leq D_0+1$ and $\psi''=2(\tilde u+S_0+M_0+1)^{-3}$,  by (\ref{2dGrad2}), we have
	\begin{align*}
		0\le&- \frac{2\epsilon_0^2 |\nabla  u|^2}{( \tilde u+S_0+M_0+1)^{3}}
		+(D_0+1) (|u_t|_{C^0}+\max_{M}\csc^2\theta_{\omega}(\chi_{u_0}))\\
		&+(D_0+1) \max_{M}|\csc\theta_{\omega}(\chi_{u_0})||\nabla \tilde u||\nabla u|^{-1}+C.
	\end{align*}
	From the above inequality, by \eqref{2dgradientassumption}, we have
	\begin{align}\label{2d0607G1}
		|\nabla  u|^2\le& C_1( 2M_0+S_0+1-S)^{3}.
	\end{align}
	By Lemma \ref{2dSestimate}, we obtain
	\begin{align}
		G(x_0, t_0)=&	|\nabla  u|^2(x_0, t_0)e^{\psi(\tilde u(x_0, t_0))}\notag\\
		\le& C_1( 2M_0+S_0+1+|S|(x_0, t_0))^{3}e^{D_0S}\le M_{1,2}.
	\end{align}
	\textbf{Case 2:}\   $\sum\limits_{i=1}^2F^{i\bar i}| \tilde u_i|^2\le \epsilon_0^2 |\nabla  u|^2$.
	
	In this case, since $\psi''>0$, by inequality (\ref{2dGrad2}), we have
	\begin{align}\label{2d060903}
		0\leq -\psi'(-u_t+F^{i\bar i}u_{i\bar i})+C\max_{M}|\csc\theta_{\omega}(\chi_{u_0})|(-\psi')\epsilon_0+C.
	\end{align}
	On the other hand, since $F^{1\bar 1}\le F^{2\bar 2}$, we have
	\begin{equation*}
		\epsilon_0^2 | \nabla  u|^2\ge  F^{1\bar 1}|\nabla  \tilde u|^2=\frac{1+\lambda_2^2}{(\lambda_1+\lambda_2)^2}|\nabla \tilde u|^2\ge \frac{1}{4\lambda_1^2}|\nabla \tilde u|^2.
	\end{equation*}
	From this inequality and \eqref{2dgradientassumption}, we get
	\begin{align*}
		\lambda_1\ge \frac{1}{4}\epsilon_0^{-1}=K_0.
	\end{align*}
	Then we can apply our Lemma \ref{2dsubsolutioninequality} to get
	\begin{align*}
		-u_t+F^{i\bar j} (u_{i\bar j}-S_{i\bar j})\le -c_0.
	\end{align*}
	Inserting the above inequality into  \eqref{2d060903}, we get
	\begin{align}\label{2d06051}
		0\le& \psi'c_0+\epsilon_0C\max_{M}|\csc\theta_{\omega}(\chi_{u_0})|(-\psi')+C\notag\\
		\le& D_0\Bigl(-c_0+\epsilon_0C\max_{M}|\csc\theta_{\omega}(\chi_{u_0})|\Bigl)+C.\notag\\
	\end{align}
	Since $\epsilon_0C\max_{M}|\csc\theta_{\omega}(\chi_{u_0})|\le \frac{c_0}{2}$, if we we choose $D_0=2c_0^{-1}(C+1)$, we get the following contradiction
	\begin{align}
		0\le -D_0\frac{c_0}{2}+C=1.
	\end{align}
	Thus this case can not occur.
	
	In conclusion, for any $(x, t)\in M\setminus\cup_i E_i$, we have $G(x, t)\le G(x_0, t_0)\le M_{1,0}+M_{1,1}+M_{1,2}$ and then we obtain the desired estimate
	\begin{align}
		|Du|^{2}(x, t)\le M_1^2 e^{D_0S(x)}=M_1\prod_{i}|s_i|_{h_i}^{-2D_0a_i}(x).
	\end{align}
\end{proof}

\begin{proposition}\label{2dsecondorder}
	There exist  uniform constant $D_1$ and $M_2$ such that for any $(x,t)\in {M\setminus\cup_i E_i\times[0, \infty)}$
	\begin{align}
		|\p\bar \p
		u|_{\omega}(x,t)\le M_2{\prod_{i}|s_i|^{-2D_1a_i}_{h_i}(x,t)}.
	\end{align}
\end{proposition}

\begin{proof}
	We consider
	\begin{align*}
		\tilde H(x, t,\xi(x))=\log(w_{i\bar j}\xi^{i}\bar{\xi}^{j})+\psi( \tilde u)
	\end{align*}
	where $\psi(\tilde u)=-D_1 \tilde u+(\tilde u +M_0+S_0+1)^{-1}$ and $\tilde u=u-S$.
	Recall $M_0$ is the uniform bound of $|u|$ in Lemma \ref{0606zoero} and $S_0$ is the upper bound of $S$. Hence we have
	\begin{align}\label{2d0606psi}
		D_1\le- \psi'\le  D_1+1 \ \ \text{and} \quad
		\psi''=2(\tilde u +M_0+S_0+1)^{-3}.
	\end{align}
	For any $T_0\in (0, \infty)$, suppose the function $\tilde H$ which is upper semi-continuous attains its maximum on $M\times [0,T_0]$ at $(x_0, t_0)$ along the direction $\xi_0=\xi(x_0)$. Since $\tilde H=-\infty$ on $\cup_i E_i$, we have $x_0\in M\setminus\cup_i E_i$. If $t_0=0$, the estimate holds since $S$ is upper bounded. Hence in the following we assume  $t_0>0$.
	
	Take holomorphic  coordinates near  $x_0$ such that \eqref{3coordinate} holds.
	Then the function\ $H(x, t)=\tilde H(x, t,\tilde\xi_0(x))$  attains its maximum on $M\times [0,T_0]$ at \ $(x_0, t_0)$.
	
	At $(x_0, t_0)$, we have
	\begin{align}
		0\le& H_t=\frac{u_{t1\bar 1}}{w_{1\bar 1}}+\psi'u_t,\notag\\
		0=&H_i=\frac{w_{1\bar 1, i}}{w_{1\bar 1}}+\psi' \tilde u_i,\label{060802}
	\end{align}
	and
	\begin{align}\label{2dPH1}
		0\le& H_t-F^{i\bar i}H_{i\bar i}\notag\\
		=&\lambda_1^{-1}(u_{t1\bar 1}-F^{i\bar i}w_{1\bar 1, i\bar i}) +\lambda_1^{-2}F^{i\bar i}|w_{1\bar 1, i}|^2 \quad \textup{(denoted by (I))}\notag\\
		&-\psi''F^{i\bar i}| \tilde u_i|^2+\psi'(u_t-F^{i\bar i} \tilde u_{i\bar i}).
	\end{align}
	By the same argument as that in section 4, $(I)$ has the following estimate
	\begin{align}
		\textup{(I)}\le&\psi'^2F^{1\bar 1}|\nabla  \tilde u|^2+ \psi'^2(1+A_1)\lambda_1^{-1}F^{2\bar 2}| \tilde u_2|^2+C.\label{2dPFI}
	\end{align}
	
	Inserting  \eqref{2dPFI} into \eqref{2dPH1}, by \eqref{2d0606psi}, we have
	\begin{align}
		0 \le& \bigl(-\psi''+(1+A_1)\psi'^2\lambda_1^{-1} \bigl)F^{2\bar 2}
		| \tilde u_2|^2+\psi'^2F^{1\bar 1}|\nabla \tilde  u|^2\notag\\
		&+\psi'(u_t-F^{i\bar i} \tilde u_{i\bar i})+C\notag\\
		\le&\bigl(-2(-S+2M_0+S_0+1)^{-3}+(1+A_1)(D_1+ 1)^2\lambda_1^{-1}\bigl)F^{2\bar 2}
		| u_2|^2\notag\\
		&+(D_1+ 1)^2|\nabla \tilde  u|^2\frac{1+\lambda_2^2}{(\lambda_1+\lambda_2)^2}
		+\psi'(u_t-F^{i\bar i}\tilde u_{i\bar i})+C.\label{PH2}
	\end{align}
	The first term is negative if we assume
	\begin{align}\label{2d0606second1}
		\lambda_1>(1+A_1)(D_1+ 1)^2(-S+2M_0+S_0+1)^3.
	\end{align}
	We further  assume
	\begin{align}\label{2d0606second2}
		\lambda_1>2K_0.
	\end{align}
	Then by Lemma \ref{2dsubsolutioninequality} and \eqref{2d0606psi}  , we have
	\begin{equation*}
		\psi'(u_t-F^{i\bar i} u_{i\bar i})\le -c_0 D_1.
	\end{equation*}
	Hence  \eqref{PH2} becomes
	\begin{align*}
		0\le 	(D_1+ 1)^2|\nabla \tilde  u|^2\frac{1+A_1^2}{(\lambda_1-A_1)^2}
		-c_0D_1+C
	\end{align*}
	or
	\begin{align*}
		\bigl(c_0D_1
		-C\bigr)(\lambda_1-A_1)^2\le (D_1+ 1)^2(1+A_1^2)|\nabla \tilde  u|^2.
	\end{align*}
	We choose  $D_1$
	\begin{align}\label{2dD11}
		D_1>{c_0^{-1}(C+1)}.
	\end{align}
	Then we have
	\begin{align}\label{2d0606second3}
		\lambda_1\le&( D_1+ 1)(1+A_1^2)|\nabla \tilde  u|+A_1\notag\\
		\le& ( D_1+ 1)(1+A_1^2)(|\nabla u|+|\nabla S|)+A_1\notag\\
		\le& ( D_1+ 1)(1+A_1^2)(M_1\prod_{i}|s_i|_{h_i}^{-M_1}+|\nabla S|)+A_1,
	\end{align}
	where in the last inequality we use \eqref{2dgradient1}.
	
	By \eqref{2d0606second1}, \eqref{2d0606second2} and  \eqref{2d0606second3},
	we obtain
	\begin{align*}
		\lambda_1\le 2K_0+(1+A_1)(D_1+ 1)^2(-S+2M_0+S_0+1)^3\\
		+ ( D_1+ 1)(1+A_1^2)(M_1\prod_{i}|s_i|_{h_i}^{-M_1}+|\nabla S|)+A_1.
	\end{align*}
	Hence
	we have at $(x_0, t_0)$,
	\begin{align}
		\lambda_1e^{\psi(\tilde u)}\le C\lambda_1e^{D_1S}\Bigr(&M_1\prod_{i}|s_i|_{h_i}^{-M_1}+|\nabla S|\notag\\
		&+(-S+2M_0+S_0+1)^3\Bigr)+C.\notag
	\end{align}
	If we choose $D_1>(M_1+1)\min_{ i}\{a_i^{-1}\}$, the above inequality has an uniform upper bound and thus we obtain the estimate \eqref{2dsecondorder}.
\end{proof}

\begin{proposition}\label{2dhighorder}
	For any compact set $K\subset M\setminus\cup_i E_i $ and positive integer $k$, there exists a uniform constant $C_{k,K}$ such that
	\begin{align}
		|u|_{C^k(K)}\le C_{k,K}.
	\end{align}
\end{proposition}

\begin{proof}
	By the complex Hessian estimate in Proposition \ref{2dsecondorder}, the  flow is uniformly parabolic. Since $\cot\theta_{\omega}(\chi_u)$ is concave, by the  Evans-Krylov theory \cite{ev1982cpam,kr1983}, we obtain the higher order estimates in $K$.
\end{proof}
As an application of Proposition \ref{2dhighorder}, we first show
\begin{proposition}\label{2duniformut}
	For any compact set $K\subset M\setminus\cup_i E_i $, $\frac{\partial u}{\partial t}$ uniformly  converges to $0$  in $K$ as $t$ tends to $\infty$.
\end{proposition}

\begin{proof}
	We first prove that  $\frac{\partial u}{\partial t}$ pointwisely  converges to $0$\\ in $M\setminus\cup_i E_i$.
	Since \begin{align}
		\mathrm{Re(CY_{\mathbb{C}}}(u(t)))&-\mathrm{Re(CY_{\mathbb{C}}}(u(0)))\notag\\
		=&
		\int_{0}^t \int_M\Bigl(\frac{\partial u}{\partial s}\Bigr)^2\mathrm{Im}(\chi_{u(s)}+\sqrt{-1}\omega)^2ds,
	\end{align}
	by Corollary \ref{2duniformboundReCY} we have
	\begin{align*}
		\int_{0}^{\infty} \int_M\Bigl(\frac{\partial u}{\partial t}\Bigr)^2\mathrm{Im}(\chi_u+\sqrt{-1}\omega)^2 dt\le C.
	\end{align*}
	Since along the flow $\mathrm{Im}(\chi_u+\sqrt{-1}\omega)^2\ge c_0\omega^2>0$, the above inequality gives
	\begin{align}\label{2duniformboundut}
		\int_{0}^{\infty} \int_M\Bigl(\frac{\partial u}{\partial t}\Bigr)^2\omega^2 dt\le c_0^{-1}C.
	\end{align}
	
	If there exists $x_0\in K$ such that $\lim\limits_{t\rightarrow \infty}\frac{\partial u}{\partial t}(x_0, t)\neq 0$, then there exists $\epsilon_0>0$ and  a sequence $\{t_i\}$ which tends to $\infty$ such that
	\begin{align}\label{2dutconverg}
		\Bigl |\frac{\partial u}{\partial t}(x_0, t_i)\Bigr|\ge \epsilon_0.
	\end{align}
	Let $U$ be a small neighborhood  of $x$ such that $U\subset M\setminus\cup_i E_i$. Then by Proposition \ref{2dhighorder}, $\frac{\partial u}{\partial t}$
	and its time and space derivative are uniformly bounded in $U\times [0, \infty)$ and thus by \eqref{2dutconverg}, there exist a small neighborhood $V\subset U$ of $x_0$ and a uniform constant $\delta>0$ such that
	\begin{align*}
		\Bigl|\frac{\partial u}{\partial t}\Bigr|\ge \frac{\epsilon_0}{2} \ \text{for any}\ (x, t)\in V\times [t_i, t_i+\delta].
	\end{align*}
	This implies
	\begin{align*}
		\int_{0}^{\infty} \int_M\Bigl(\frac{\partial u}{\partial t}\Bigr)^2\omega^2dt\ge& \sum_{i=1}^{\infty}\int_{t_i}^{t_i+\delta}\int_V
		\Bigl(\frac{\partial u}{\partial t}\Bigr)^2\omega^2dt\\
		\ge& \sum_{i=1}^{\infty}\delta\frac{\epsilon_0^2}{4}\mathrm{vol_{\omega}} (V)=\infty,
	\end{align*}
	which contradicts with \eqref {2duniformboundut}.
	Hence $\frac{\partial u}{\partial t}$ point-wisely  converges to $0$ in $M\setminus\cup_i E_i$.
	
	Let $K\subset \cup_{j=1}^N B_{r}(x_j)\subset M\setminus \cup_i E_i$. We can apply the differential Harnack inequality for $\frac{\partial u}{\partial t}$ in every $B_{r}(x_j)$ to  prove that   $\frac{\partial u}{\partial t}$ converges in any compact subset $K$ uniformly to $0$. 
\end{proof}

\subsection{Proof of Theorem \ref{2dtheorem1}.}
Similarly as the proof by Fang-Lai-Song-Weinkove \cite{flsw2014apde} and Takahashi \cite{ta2021cvpde}, we have

\begin{lemma}\label{2dfunctional}
	Let $\{u_i\}$ be a sequence of smooth functions satisfying $\chi_{u_i}-\cot B_1 \omega>0$ and $|u_i|_{C^0}\le C$ for $C>0$. Let $u^{\infty}$ be a bounded $(\chi-\cot B_1 \omega)$-PSH function on $M$. Let $Y$ be a proper subvariety of $M$.  Assume that  $u_i$ converges to $u^{\infty}$ in $C_{loc}^{\infty}(M\setminus Y)$ as $j\rightarrow \infty$. Then $CY_{\mathbb{C}}(u^{\infty})$ and $\mathcal{J}(u^{\infty})$ are well-defined. Moreover,
	\begin{align*}
		\lim_{i\rightarrow \infty} \mathrm{Im}(\textup{CY}_{\mathbb{C}}(u_i))=&\mathrm{Im}(\textup{CY}_{\mathbb{C}}(u^{\infty})),\\
		\lim_{i\rightarrow \infty} \mathrm{Re}(\textup{CY}_{\mathbb{C}}(u_i))=&\mathrm{Re}(\textup{CY}_{\mathbb{C}}(u^{\infty})),\\
		\lim_{i\rightarrow \infty}\mathcal{J}(u_i) =& \mathcal{J}(u^{\infty}).
	\end{align*}
\end{lemma}

\begin{proof}[Proof of Theorem \ref{2dtheorem1}]
	By the $C^0$ estimate proved in Proposition \ref{2dzeroorder}, there exists a sequence $\{t_i\}$ such that $u(\ , t_i)$ converges to a function $u^{\infty}\in L^{\infty}{(M)}$. By the $C^k$ estimates in Proposition \ref{2dhighorder}, by passing a subsequence (for convenience we still denote by ${t_i}$), $u(\ , t_i)$ smoothly converges to $u^{\infty}$ in any compact subset of $M\setminus\cup_i E_i$ and thus $u^{\infty}\in C^{\infty}{(M\setminus\cup_i E_i)}$.
	Since $\chi_u>\cot {B_1}\omega$, then $\chi_{u^{\infty}}-\cot {B_1}\omega$
	is a K\"{a}hler current and is smooth in $M\setminus\cup_i E_i$. By Lemma \ref{2dfunctional} and Lemma \ref{IMCY},
	we have $ \mathrm{Im}(\textup{CY}_{\mathbb{C}}(u^{\infty}))=\mathrm{Im}(\textup{CY}_{\mathbb{C}}( {u_0}))$.
	
	By Proposition \ref{2duniformut}, $u^{\infty}$ satisfies \eqref{2d} in $M\setminus\cup_i E_i$ and
	then $\theta_{\omega}(\chi_{u^{\infty}})=\theta_0$ on $M\setminus\cup_i E_i$.
	We can define $\chi_{u^{\infty}}^2$ and $\chi_{u^{\infty}}\wedge \omega$ as finite measures on $M$ such that they do not charge pluripolar subsets.  Thus ${(\chi_{u^{\infty}}+\sqrt{-1}\omega)^2}$ is well-defined and  $u^{\infty}$ satisfies the equation \eqref{2d} on $M$ in the sense of currents. Moreover, ${u^{\infty}}$ is  $\tilde \chi$-PSH on $M$ and satisfies the equation  \eqref{2dmongeampere2} in the sense of currents.
	
	Finally, by the $C^{\infty}_{loc}(M\setminus\cup_i E_i)$ uniform estimate of $u(t)$ and the uniqueness of the equation  \eqref{2dmongeampere2}, similar as the argument in \cite{flsw2014apde}, we have $u(t)$ converges smoothly to $u^{\infty}$ on $M\setminus\cup_i E_i$.
\end{proof}

\subsection{$\mathcal J$-functional.}
As an application of our flow, we  prove the lower bound of the $\mathcal{J}$-functional in the following  set.
\begin{align*}
	\mathcal{H}_{B_1}=\left\{w\in C^{\infty}(M, \mathbb{R}): \theta_{\omega}(\chi_w)\in (0, B_1)\right\}.
\end{align*}

\begin{corollary}\label{2dcorollary}
	Let $(M, \omega)$ be a compact K\"{a}hler surface and $\chi$  a closed real $(1,1)$ form.
	Assume that $\theta_0\in(0,\pi)$ and $\chi \ge\cot\theta_0\omega $,
	the $\mathcal{J}$-functional  is bounded  from below in $\mathcal{H}_{B_1}$ for any $B_1\in(\theta_0, \pi)$.
\end{corollary}

\begin{proof}
	For $u_0\in \mathcal{H}_{B_1}$, let $u(t)$ be the solution of our flow $u_t=\cot\theta_{\omega}(\chi_u)-\cot\theta_0$ with $u(0)=u_0$. By Theorem \ref{2dtheorem1}, $u(t)$ converges to a bounded function $u^{\infty}$ solving \eqref{2dmongeampere2}. Since $\mathcal{J}$ is decreasing along the flow, we have
	\begin{align*}
		\mathcal{J}(u_0)\ge \lim\limits_{t\rightarrow \infty}\mathcal{J}(u(t))=\mathcal{J}(u^{\infty}).
	\end{align*}
	Let $v$ be a weak solution of \eqref{2dmongeampere2} in Lemma \ref{2dmongeampere1}. By the uniqueness, there exists a constant $c_0$   such that $u^{\infty}=v+c_0$. Since $\mathcal{J}(u^{\infty})=\mathcal{J}(v)$, we get
	\begin{align*}
		\mathcal{J}(u_0)\ge \mathcal{J}(v).
	\end{align*}
\end{proof}

{\bf Acknowledgements.} Zhang would like to thank Prof. Xi-Nan Ma for constant help and encouragement. Fu is supported by NSFC grant  No. 12141104 and 11871016. Zhang is supported by NSFC grant  No. 11901102.

\end{document}